\documentclass[reqno,centertags, 12pt]{amsart}
\usepackage[alphabetic,msc-links]{amsrefs}
\usepackage{amsfonts,amssymb,amsmath,amsthm,graphicx,overpic,color}
\usepackage{amsthm,amscd,amssymb,mathtools}
\usepackage[width = 6 in]{geometry}
\usepackage{float}
\usepackage{lineno}
\usepackage{wrapfig}
\usepackage[labelformat=empty]{caption}

\DeclareMathAlphabet\gothic{U}{euf}{m}{n}

\newcommand{\bdone}{{\boldsymbol{1}}}

\begin{document}
\title[Remembrances of Derek William Robinson]{Remembrances of Derek William Robinson\\
	25 June 1935 -- 31 August 2021}

\author[Editors: Michael Barnsley, Bruno Nachtergaele, Barry Simon]{Editors: Michael Barnsley, Bruno Nachtergaele, Barry Simon\\ Contributors: Alain Connes, David Evans, Giovanni Gallavotti,\\ Sheldon Glashow, Arthur Jaffe, Palle Jorgensen, Aki Kishimoto,\\ Elliott Lieb, Heide Narnhofer, David Ruelle,\\ Mary Beth Ruskai, Adam Sikora, A. F. M. ter Elst}
\date{}
\maketitle

\textbf{Editor's Notes}: This arXiv submission is posted primarily to make available some additional material to what we prepared for the Derek Memorial artcle published in the Notices of the AMS:

\emph{Remembrances of Derek William Robinson}, June 25, 1935–August 31, 2021, Notices A.M.S  70 (2023), 1252-1267.

We start (with the permission of the AMS) with the text of our submission close to the final published version (which you can get at \\
https://www.ams.org/journals/notices/202308/noti2765/noti2765.html?adat\\=September\%202023\&trk=2765\&galt=none\&cat=commentary\\ \&pdfissue=202308\&pdffile=rnoti-p1252.pdf).  After that we include some autobiographical notes (page \pageref{D1}) that Derek dictated to Louisa Barnsley (whom we thank) during his final illness.  Finally some unpublished notes (page \pageref{D2}) Derek wrote on the history of his seminal monographs with Bratelli.  We thank Marion Robinson for permission and encouragement to post these items.

\section*{Introduction by the Editors}

Derek William Robinson (June 25, 1935 – August 31, 2021) was a founding figure in the modern mathematical approach to the foundations of statistical mechanics together with Roland Dobrushin, Rudolph Haag, Oscar Lanford and especially David Ruelle. With Ola Bratteli, he wrote what has long been the standard monograph on the use of operator algebra methods in quantum statistical mechanics~\cite{BratRob} and an important resource for those studying $C^*$-algebras, an area in which he likewise made significant contributions.  Later he worked on analysis on Lie Groups where he also wrote two monographs~\cite{RobLie, DER03}.

He was born in Southern England; his father was from rural Northern Ireland.  He entered Oxford in 1954 and graduated with honors in mathematics in 1957.  Because the head tutor in Maths at Wadham College was an applied mathematician, Derek wound up choosing hydrodynamics and quantum mechanics as his special subjects.  He was offered a graduate scholarship in Physics at Oxford even though as he commented {\em the only openings they had in theoretical physics were in nuclear physics, down and dirty stuff, so I ended up having to work through this stuff in which I had absolutely no interest}.

\begin{figure}[h!]
 \centering \includegraphics[width=0.4\textwidth]{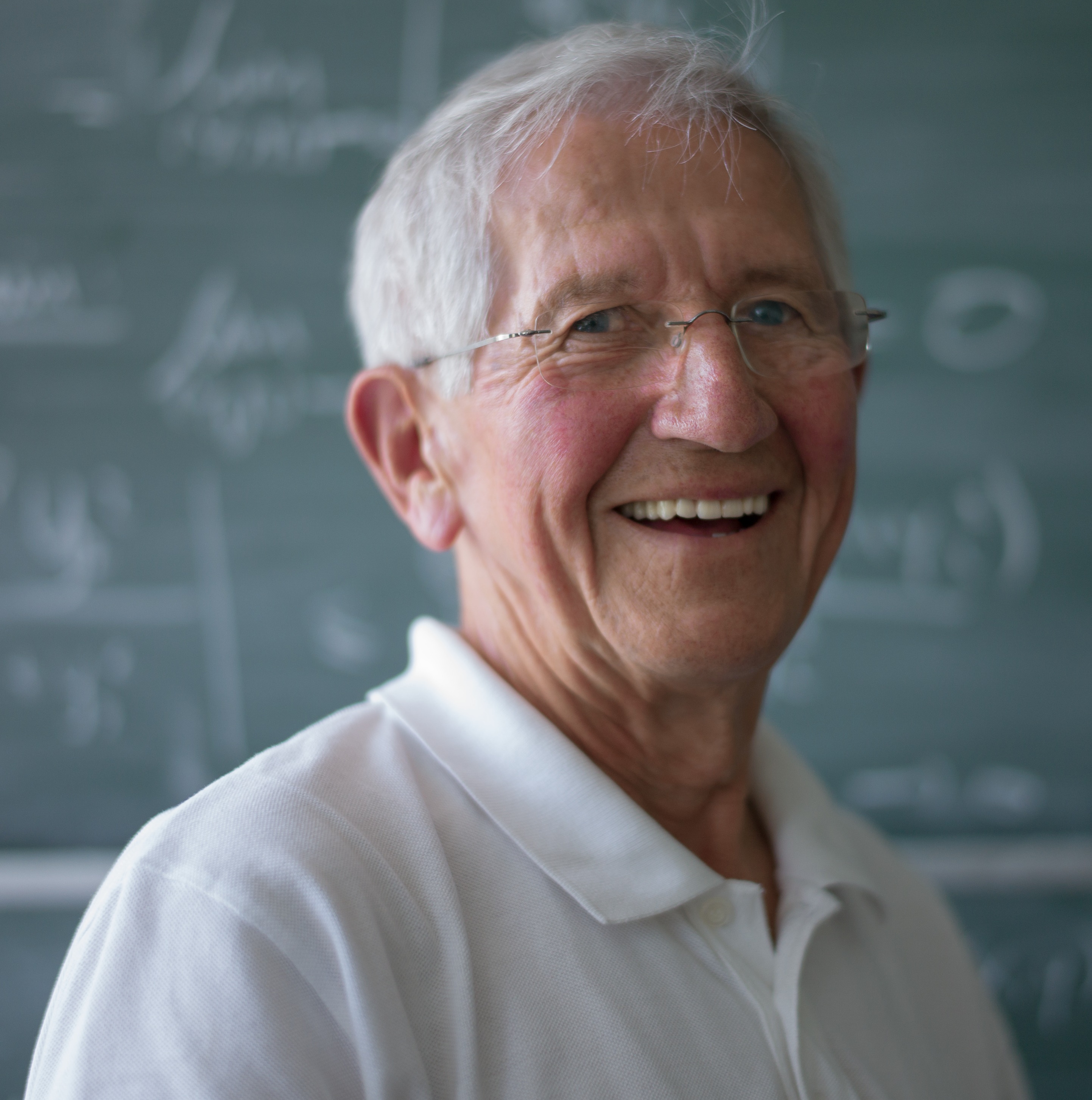}
   \caption{Derek Robinson (1935-2021)}
\end{figure}

So, although his PhD. thesis was in nuclear theory, he mainly studied quantum field theory and particle physics.  His directions were heavily influenced by summer schools he attended in those subjects in Naples and Edinburgh where he met and became good friends with two future Nobel Laureates, Tini Veltman and Shelly Glashow.  Indeed, Veltman was best man at Derek’s wedding and the Glashows spent part of their honeymoon in the Robinson’s home in the south of France.

One of the lecturers in the school in Naples was Res Jost who was Pauli’s successor at the ETH in Zurich.  Derek was so struck by his research that he took a NATO postdoc in Jost’s group in 1960-62 immediately after he finished his PhD.  He was allowed to take this fellowship at ETH even though Switzerland is not a NATO country!  One of the other postdocs in Jost’s group was Ruelle and a post-doctoral visitor was Araki who lectured about von Neumann algebras.

Derek’s second postdoctoral mentor was Rudolph Haag who together with Kastler invented the $C^*$-algebraic approach to quantum field theory and to quantum statistical mechanics.  Derek worked with Haag for roughly four years (1962-1966), two years at the University of Illinois where Haag was a professor, then 18 months in Munich while Haag was on leave.  When Haag returned to Illinois at the end of 1965, Derek stayed in Europe spending 3 months at IHES at Ruelle’s invitation and 6 months in Marseille as Kastler’s visitor.

Robinson then moved to a position within the theory group at CERN during 1966-68.  While the group mainly focused on high energy physics, during this period, Derek mainly did statistical mechanics (as he had during the visits to IHES and Marseille) doing seminal work we’ll describe below.  This work led to several tenure offers.  While in Geneva, Derek met his wife Marion, the daughter of a diplomat from New Zealand who had spent her youth in Canada and mainly Geneva.

\begin{figure}[h!]
 \centering \includegraphics[width=0.4\textwidth]{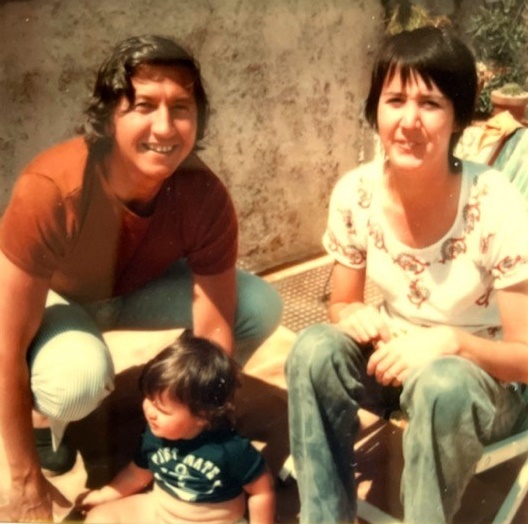}
   \caption{Derek and Marion in Bandol}
\end{figure}

Derek and Marion decided to accept the offer Kastler had arranged for a Professorship in Theoretical Physics from the University of Aix-Marseille.  Derek was one of the first non-native French professors because one of the results of the 1968 student uprising was to make it possible for non-French people to have permanent public service jobs including permanent university appointments!  Derek and Marion built a house in Bandol, a town on the Mediterranean coast east of Marseille where the Kastlers lived.

\begin{figure}[h!]
 \centering \includegraphics[width=0.4\textwidth]{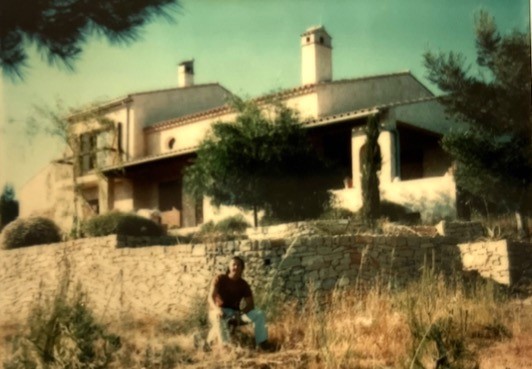}
   \caption{The House in Bandol}
\end{figure}

During the time in Marseille, Derek’s research shifted more towards problems in the theory of operator algebras and he began working with Ola Bratteli with whom he’d work for over 40 years including their book which they initially planned to be 300-400 pages but turned out to be two volumes totaling more than 1050 pages!

In early 1976, Angas Hurst (the H of the GHS inequalities) lured Derek to Adelaide to lecture in the “Summer” Research school.  Derek was taken with Australia and he and Marion began thinking about a permanent move.  He applied for a Professorship in Pure Mathematics at the University of New South Wales where he was from 1978 until 1982.  Derek then moved to a Professorship in the Institute for Advanced Studies at the Australian National University in  Canberra, a position he held for 19 years until his retirement. While he became emeritus in 2000, he continued to have a research grant and wrote papers until shortly before his death.

Derek was a serious competitive cyclist and at one point wrote that {\em you don't necessarily peak as a cyclist when you are young, I found I was faster at time trials at
the age of 50 than I had been as an eighteen-year-old. In 2002 I was very pleased to win the World Masters Games 20k cycling race in my age range. This was a time trial and I was the last person to start as I was the top seed and I completed the course in a little over 32 minutes. In some ways I consider this my greatest achievement!}

\begin{figure}[h!]
 \centering \includegraphics[width=0.3\textwidth]{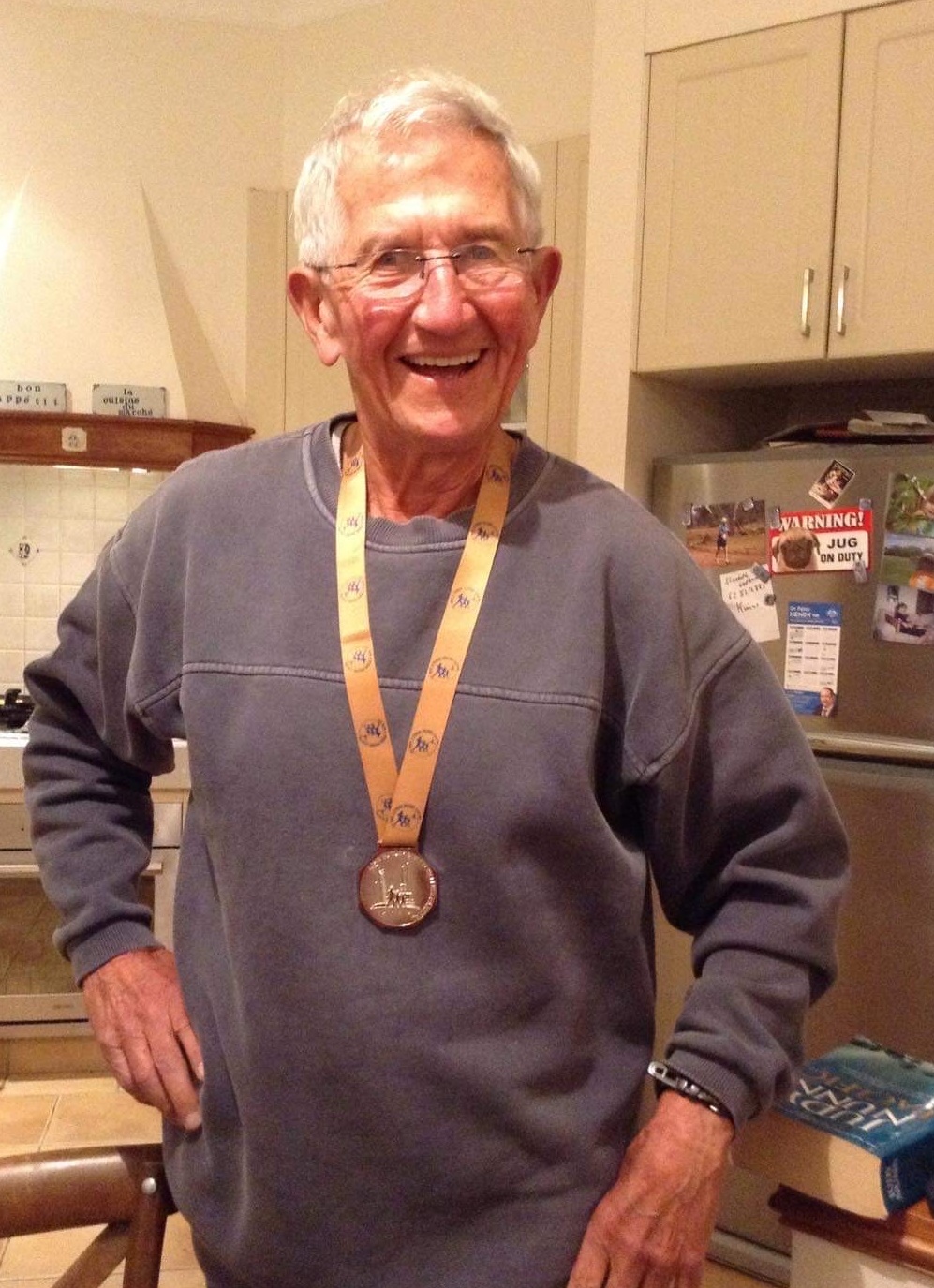}
   \caption{Derek with his cycling medal, 2002}
\end{figure}

We turn to a description of some of Derek’s best known research.  In the mid-1960's there was a major paradigm shift in which Derek was a major figure so that the study of systems in infinite volume took over and new mathematical realizations, problems and objects came to the fore.  There was considerable resistance in the theoretical physics community even though now, the picture is so accepted that it is hard to realize there was ever an issue.  Derek's work focused on issues special to the quantum case although his paper with Ruelle on entropy \cite{RobRuEnt} in the classical setting was a factor in the classical framework which Derek sought to extend to quantum systems.

One of his key discoveries concerned the invention of asymptotically abelian actions.  For classical systems, an important property of infinite volume equilibrium states is the uniqueness of decomposition into pure phases - in mathematical terms, the crucial underlying fact is that the translation invariant probability measures are a Choquet simplex whose extreme points are the ergodic measures.  It inherits this from the fact that all probability measures on a compact Hausdorff space are a simplex.

For quantum systems, there is an issue in extending this notion: the set of all states on a $C^*$ algebra is a simplex if and only if the algebra is abelian and, of course, general quantum systems are non-abelian!  Derek discovered the idea of an asymptotically abelian action by translations – those with the property that for any pair, $A, B$ of operators, $A$ and the distant translates of $B$ asymptotically commute.  The set of states invariant under such an action are a simplex!  This idea was developed by Derek in three papers with Doplicher and Kastler \cite{AsymAbel1} and one with Ruelle.

A second foundational question is the existence of an infinite volume limit for dynamics.  While dynamics plays a minor role in classical spin systems, it is central to quantum spin systems, for example, via the KMS boundary condition.  For both cases, one is mainly interested in interactions with an upper bound on the number of bodies involved in the potentials but since one wants arbitrarily many bodies and it is technically very convenient to have a Banach space.  What Derek found, as part of a series of three papers setting the framework for quantum spin systems \cite{QSpin}, was a smaller Banach space than for the classical case where one could always define dynamics.

\begin{figure}[h!]
 \centering \includegraphics[width=0.2\textwidth]{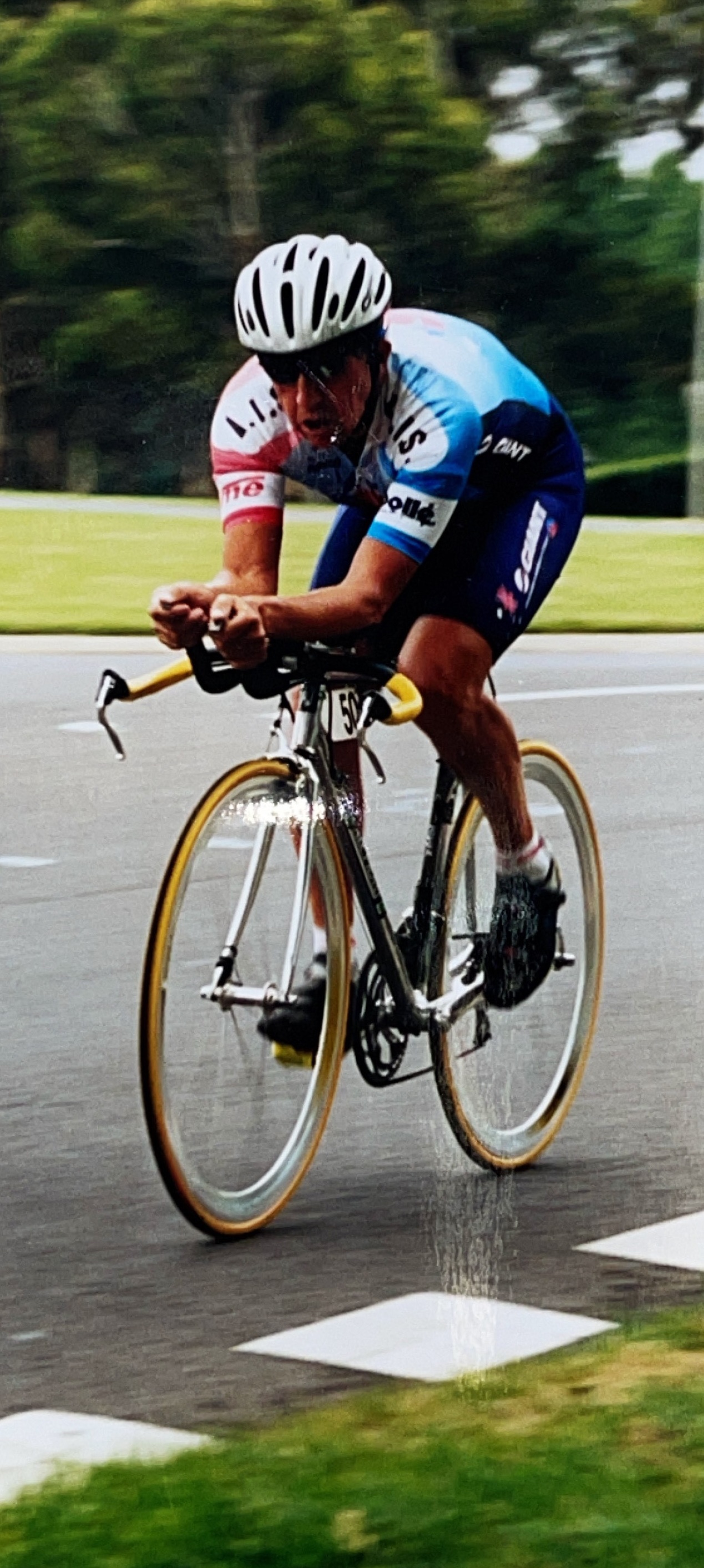}
   \caption{Derek on His Bike}
\end{figure}

A third component involved finding to what extent the entropy ideas he and Ruelle had studied for classical systems \cite{RobRuEnt} extended to quantum systems.  His classic paper with Lanford \cite{LandRob} settled many questions but also raised others including the question of strong subadditivity that we’ll hear about in the contributions below of Lieb and Ruskai.

We should mention one other item before leaving statistical mechanics: the bound on speed of propagation in quantum spin systems which Derek and Elliott Lieb found \cite{LR72} with the following remarkable statistic: during its first $30$ years, it averaged fewer than one citation per year but after Hastings \cite{Hast1} in 2004 discovered its relevance to the then fledgling field of quantum information theory, there was change, so much so that recently it has averaged about 3 citations per {\em week}.  The contribution of Nachtergaele below discusses this further.

\begin{figure}[h!]
 \centering \includegraphics[width=0.4\textwidth]{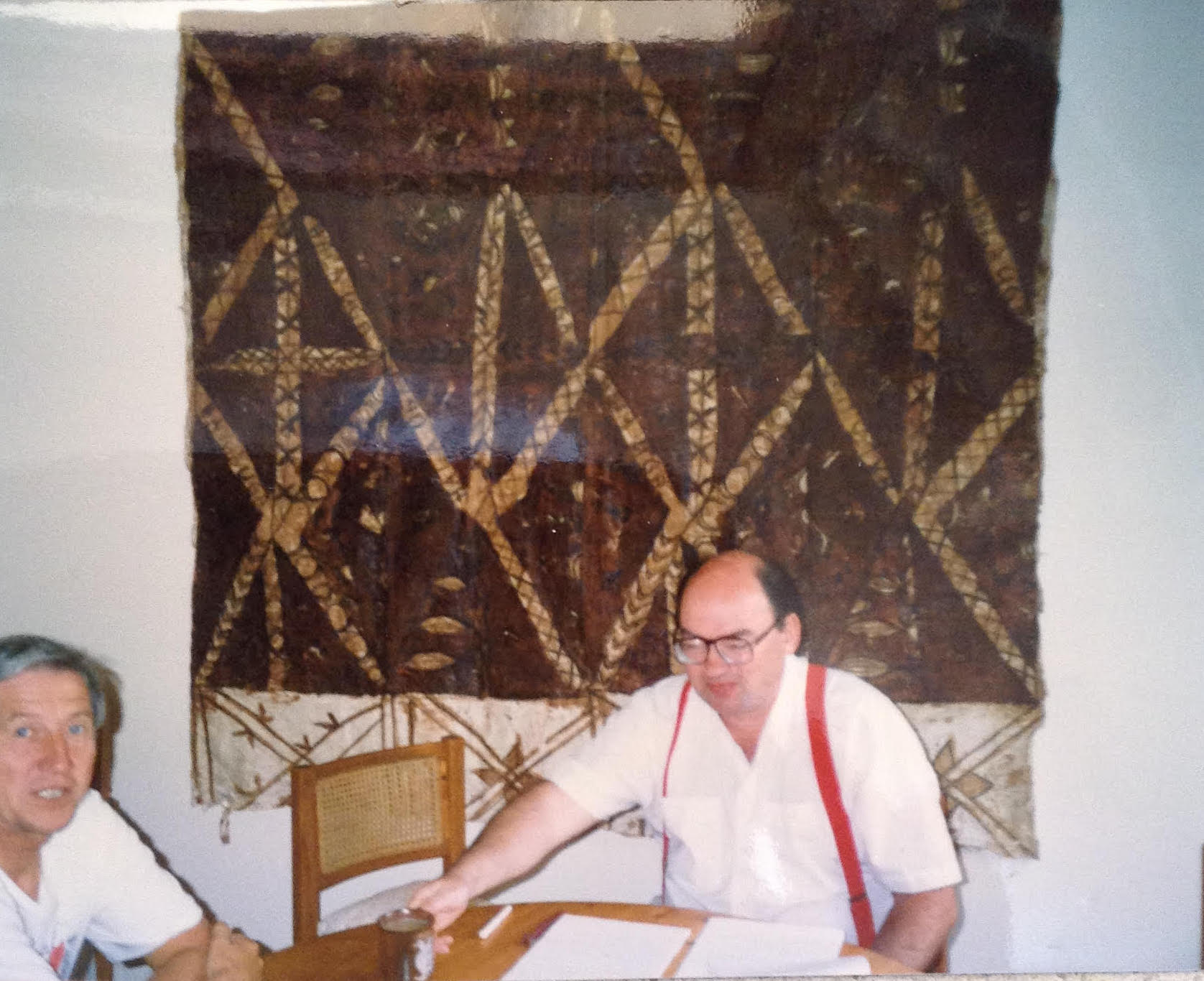}
   \caption{Derek with Ola Bratteli}
\end{figure}

If there is a leitmotif to Derek's later work in $C^*$-algebras and the theory of elliptic operators on Lie groups (discussed further by ter Elst), it is dynamics, a subject which Derek addressed for quantum spin system in \cite{QSpin}.  With Bratteli \cite{DerivCstar}, he explored the theory of unbounded derivations on operator algebras.  Derek and Ola made many striking discoveries starting with the fact that such derivations are not always closable.  Derek wrote about his collaboration with Ola Bratteli in the \emph{Notices} memorial to Bratteli\cite{RobBrat}.

For more than fifty years, Derek produced beautiful and significant mathematics.

\emph{The editors would like to thank Palle Jorgensen and Beth Ruskai for help in planning and Louisa Barnsley and Marion Robinson for information about Derek's life.  We relied on some autobiographical notes of Derek's that Louisa transcribed.  We expect to place them with a final draft of this article on the arXiv.}

\section*{Michael Barnsley}

I first noticed Derek maybe fifteen years ago. He was at the back of the tearoom of the Mathematical Sciences Institute (MSI) at Australian National University (ANU), writing neatly on the blackboard, showing mathematical things to Adam Sikora. Teas were on Thursdays before the weekly colloquium; many mathematicians, some famous, flowed in and out, leaving formulas on the boards and sometimes an empty wine bottle, marking the timelessness of mathematics and, it seems to me, of Derek. In 2014 he presented an invited paper at a conference, New Directions in Fractal Geometry; his talk concerned uniqueness of solutions to diffusion equations in media with fractal boundaries, and was illustrated by Louisa Barnsley. Over the following years he wrote a number of papers in this direction. And we learnt over cups of tea and biscuits, over the years, stopping by my office, of his days growing up in England, O-levels, A-levels, Oxford Entrance exam, and how he followed his own path.

\begin{figure}[h!]
 \centering \includegraphics[width=0.25\textwidth]{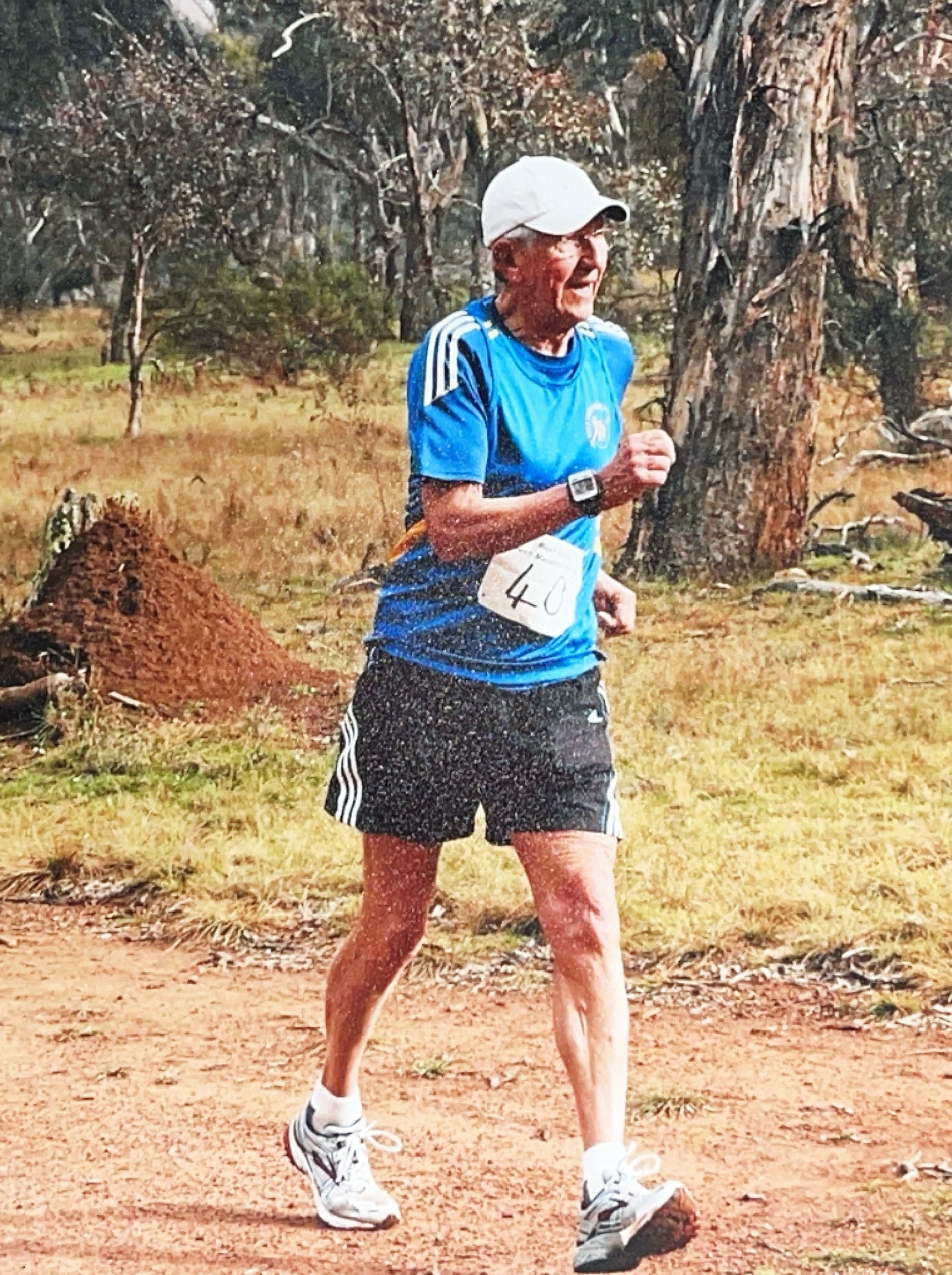}
   \caption{Derek as a  speed walker}
\end{figure}

In the years Louisa and I knew him, in addition to his steady mathematical work, Derek was proud of being the fastest speed walker in his age group, competing up and down Mt. Anslie, near ANU. But his real interest was mathematics, and he tracked the growing list of citations to his earlier papers, talked about his current work, and sometimes asked questions about fractal geometry, such as how best to understand the relationship between the inside and the outside of a Koch curve. He had answered his own question the next morning, before I could hazard a description.

At the end, we went to his daughter's beautiful apartment in Kingston, where we were greeted by Marion and Sasha with homemade cake and tea. We sat beside Derek in his special hospital bed, installed next to huge plate glass windows, looking across Lake Burley Griffin at cows in long grasses and rushes. The sunlight streaked in. Over a number of visits, sometimes pushing him in a wheelchair beside the lovely lake, he told us the story of his early days. I told him I thought the Notices would publish a memorial article about him. We organized a one-day conference entitled \textquotedblleft The Mathematical World of Derek Robinson\textquotedblright\ held on 30 July 2021. The speakers included Alain Connes, Elliott Lieb, David Ruelle, Mary Beth Ruskai, Bruno Nachtergaele, Tom ter Elst and Juha Lehrbeck, as well as local members of the MSI. Derek had never met Juha, but was thrilled to see him on Zoom at the meeting. In his penultimate paper \cite{Rob21A}, Derek included a little fractal picture.


\section*{Alain Connes}

This short text will try to give an idea of how Derek was, while in Bandol in the seventies where Daniel Kastler had gathered around his villa in Bandol a number of theoretical physicists devoted to quantum statistical mechanics and algebraic quantum field theory.

Derek Robinson was playing a key role among them by his remarkable mathematical and physical intuition, plus a conjunction of great expertise in functional analysis and deep knowledge of quantum physics.

Derek had decided to have his house in an isolated location not far from Daniel’s house but not easily accessible by car. Besides his shining scientific talent, Derek was a joyful person with lots of wit and I remember an occasion, on a dinner in his house for Thanksgiving, when at some point Derek disappeared from the dinner table, went to the Kitchen, and scared all of us by a loud scream. When asked by his shivering wife why he shouted like that, Derek answered:
{\em There is a Bird in the Oven\dots }.

At some point Derek decided to build by himself a swimming pool next to his isolated house. After building a high wall, it turned out that by accident the wall fell on him! He was alone at that point and managed, with a broken back, to crawl to a phone, call an ambulance which then took him through his bumpy road to the next hospital. After a few months in a cast Derek reappeared and his first words to me were:
{\em I just spent months in the unit ball of a Banach space}.


\section*{David Evans}

When Derek was getting settled at ANU, he invited me to visit as we had mutual interests in dynamical semigroups of positive maps on operator algebras. Consequently during a sabbatical in 1982--83, I went to Canberra for four months in late 1982 at the Research School of Mathematics which then shared a building with the Research School of Theoretical Physics. Charles Batty and Marinus Winnink were also visiting. We all met most mornings in Derek’s office discussing not only our dynamical semigroups but a variety of topics including the $C^*$-operator algebra approach to the Ising model which Winnink and I had worked on separately with John Lewis at Dublin. At morning coffee one day, Rodney Baxter brought a copy of his book {\em Exactly Solved Models in Statistical Mechanics} fresh from the press. I did not then realise what a profound impact this book would have on my future work and that it was related to a preprint I was carrying with me which Vaughan Jones gave me around the same time on his new index of subfactors. The 1982 visit was the start of numerous bilateral visits between Canberra and Warwick. This initial visit to ANU was also life changing for me in that I met Pornsawan who would become my wife.

\begin{figure}[h!]
 \centering \includegraphics[width=0.40\textwidth]{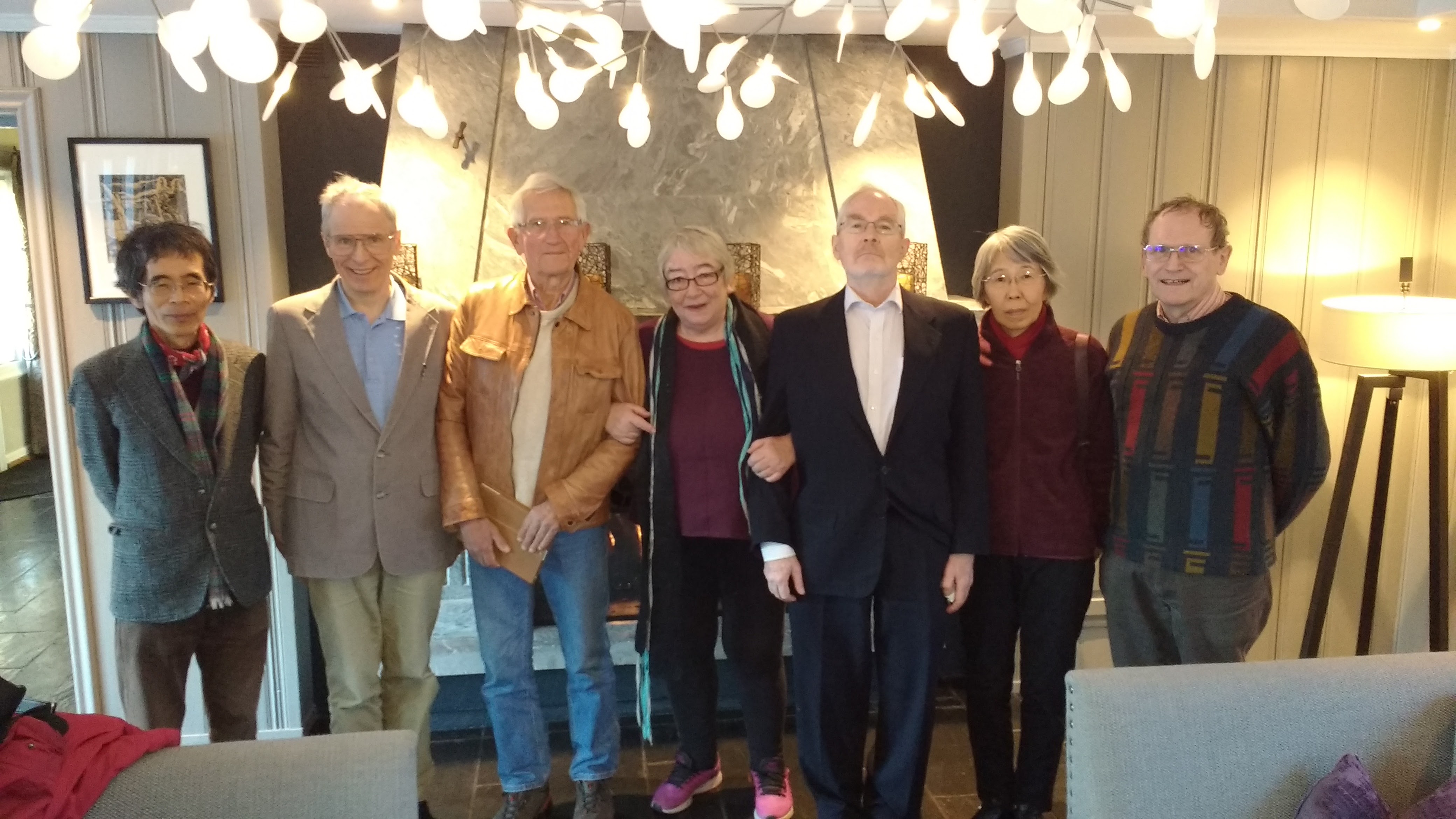}
   \caption{Lysebu, Oslo, 2017. From left Aki Kishimoto, Palle Jorgensen, Derek Robinson, Tone Bratteli (Ola's sister), George Elliott, Reiko Kishimoto, and David Evans}
\end{figure}

After Canberra, I visited Huzihiro Araki for few months in 1982-83. We wrote a paper on understanding the classical two dimensional Ising phase transition through the ground states on the one dimensional quantum system obtained from the transfer matrices and started thinking about extending this approach to the Potts model. At the end of that sabbatical year, when visiting Winnink in September 1983, I realised that the algebraic relations of the Temperley-Lieb projections within the transfer matrix formalism for the Potts model, which I found in Baxter’s book, were exactly the same as in the family of Jones projections in his tower of subfactors. Moreover the Pimsner-Popa representation of these projections, which I had just learnt about at a conference in Romania in August 1983, was the same as that discovered by Temperley and Lieb when showing an equivalence between the Potts model and an ice-type model. This and subsequent developments were triggered by Derek’s invitation to Canberra. On one visit to Warwick in 1984, Derek bought near London a high performance racing bike. On his next visit, to the 1986-87 Warwick Symposium I ran on Operator Algebras, the company arranged for a car to pick him up on arrival at Heathrow airport to bring him to the store for his next cycle purchase. Derek was an integral part of that Symposium, and he left an indelible memory not only through his scientific contributions but also when turning up at the department in lycra. Derek's friendship, sense of humour and encouragement to discuss mathematical physics on a broad spectrum are memories I treasure.


\section*{Giovanni Gallavotti}

I have been greatly influenced by Derek Robinson's works and ideas that he developed at IHES in the years 1966-1968, which in the early 1970s I soon conveyed to my students and collaborators, \cite{AsymAbel1, QSpin, RobRuEnt}.  I remain very grateful  for his kindness and patience.

My recollection of Derek goes back to 1967 when, at IHES, together with Salvador Miracle-Sol\'e we had the exiting experience of realizing the high temperature analyticity of the
Ising model thermodynamic properties, and, shortly afterwards, the same idea was applied to the Heisenberg model, making use of the path integral representation found by Jean Ginibre.

\begin{figure}[h!]
 \centering \includegraphics[width=0.25\textwidth]{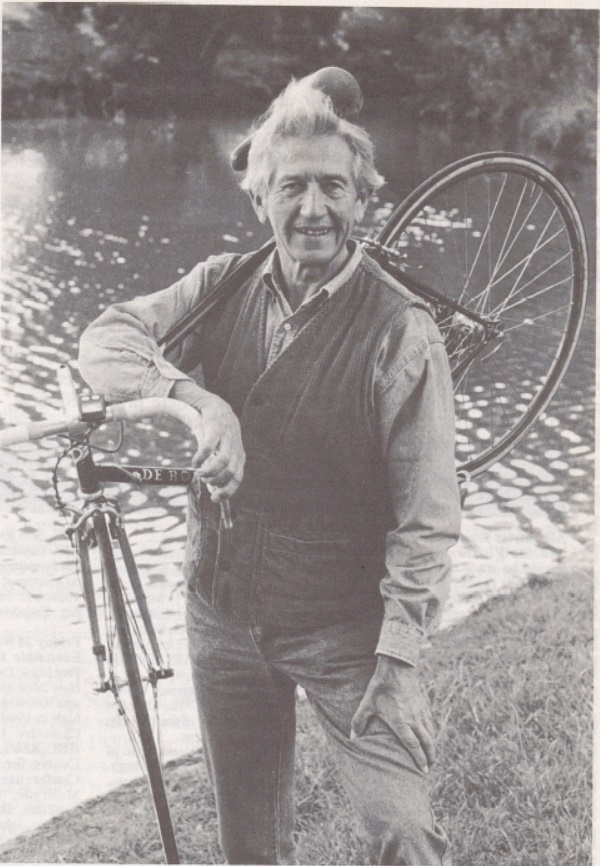}
   \caption{Derek with His Bike, 1995}
\end{figure}

Then he moved to Marseille and in 1976 I had a brief but intense relation with him: he and Ola Bratteli had spotted a serious mistake in a preprint by Mario Pulvirenti and me about the classical version
of the KMS condition, and we had to work very hard to find a correction which eventually was included in the fundamental book by Bratteli and Robinson \cite{BratRob}.

Shortly afterwards Derek moved to Australia where I met him again, briefly, in 1997 in Brisbane at the XIIth International Congress of Mathematical Physics: physically he had not changed much and he confirmed that he was still cycling as his (professional level) passion for this sport had not dwindled. I could only regret that the physical distance between Italy and Australia is too large to
allow resuming the intense collaboration of 1967 and the flash interaction of 1976.


\section*{Sheldon Glashow}

Derek was my close friend for decades. Soon after he moved to Australia, he told me his only regret was not doing so earlier. We overlapped throughout my two European sabbatical semesters, one in Geneva in the 60s, the other in Marseille in the 70s.  We never collaborated scientifically, but we thoroughly enjoyed the release of Sgt. Pepper's Lonely Hearts Club Band in 1967. We have only once or twice seen each other since his continental shift. Let me recall a couple of ancient anecdotes, possibly a bit colored by time.

Derek and Marion had been an item for some time before they chose to marry. Living in France, they found it very difficult to collect the required documents. Appealing to the British consul in Marseille, they learned they needed no further documents. Because they were Commonwealth citizens, the consul could marry them there and then.  And so he did, throwing in the Champagne.

Upon one of Derek's visits to the States in early 1972, I invited him to dinner at the home of my then girlfriend. After the gourmet dinner Joan prepared, Derek and I enjoyed every drop of the Marc de Provence he had brought me from France. He collapsed on the couch for the night, while I proposed to Joan.  We wed a few months later. This year, Joan and I celebrated our golden wedding anniversary.

We spent our honeymoon weeks partly in Paris, but mostly with Derek and Marion in the country house they built.   They were still haunted by the French bureaucracy. Just as their new kitchen was completed, the rules governing utility dimensions was changed and their new dishwasher could not be installed. Nonetheless, Joan and I had a marvelous honeymoon ensconced with our dear friends deep in the woods surrounding beautiful downtown Bandol, where the wine was cheap but good. Today, it has become expensive but excellent. Plus ca change.

Derek Robinson will always remain a part me. He was a remarkable scientist and the dearest of friends.


\section*{Arthur Jaffe}

Encouraged by my fellow student Gavin Wraith at Cambridge, UK, I travelled in the summer of 1961 as an “observer” to a school in theoretical physics held in Hercegnovi, a beautiful city on the Adriatic coast, then part of Yugoslavia. While I was too young to attend the school officially, it was serendipity to be present, for in Hercegnovi I met three students from the ETH: Derek, Klaus Hepp, and David Ruelle. Memorable teachers who interacted with everyone included Kurt Symanzik, Tulio Regge, Walter Thirring, Maurice Jacob, and Andr\'{e} Martin. I am sure that ebullient Derek got to know them all.

Derek and I overlapped again at the meeting organized by Dick Kadison at Louisiana State University. To everyone’s dismay, the Japanese mathematician Minoru Tomita arrived to announce his astounding insight, now known as modular theory. This marked the 1967 Baton Rouge meeting as one of the most influential mathematical conferences ever held in the US. I am unsure whether anyone understood Tomita’s talk at the time, but Masamichi Takesaki was fascinated by it and spent the next years as a postdoctoral fellow with Kadison in Pennsylvania—working out the details and writing his famous 1970 book.  My friend Sergio Doplicher (whom I met when we both students at the IHES during the 1963-4 special year on mathematical quantum field theory) also came to Baton Rouge.  Derek, Sergio, and I spent one day together after the meeting, mostly as tourists in New Orleans, when I took the Polaroid shown.

\begin{figure}[h!]
	\centering \includegraphics[width=0.35\textwidth]{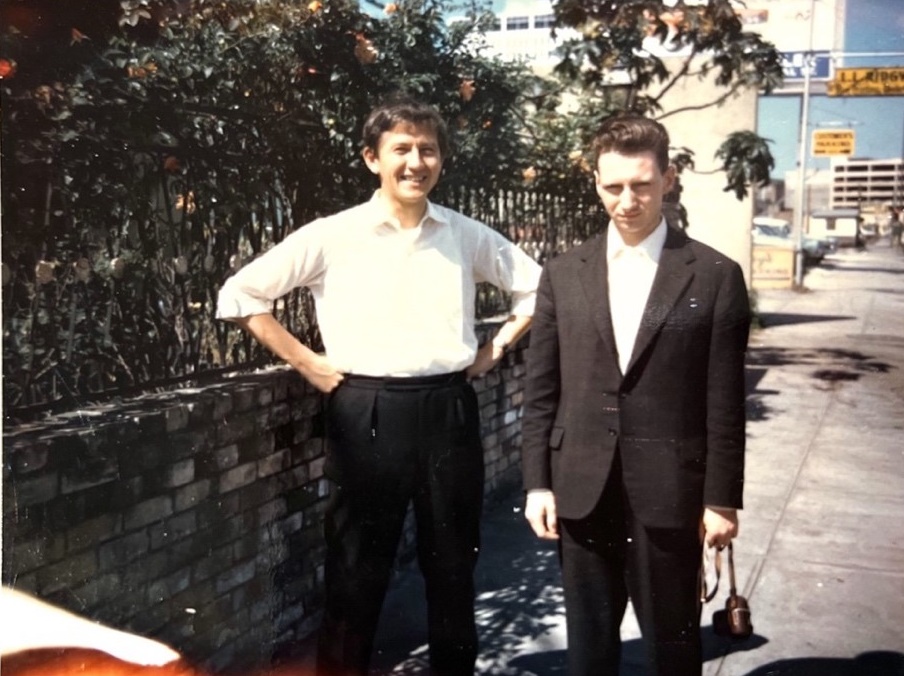}
	\caption{Derek with Sergio Doplicher, New Orleans,\\ after the 1967 Baton Rouge Conference}
\end{figure}

In 1979-81, Derek and Ola Bratteli published their monumental two-volume monograph, on \emph{Operator Algebras and Quantum Statistical Mechanics}. Here they give a beautiful derivation of the Tomita-Takesaki theory, along with explaining its importance and relation to other subjects, including its relation to the Kubo-Martin-Schwinger condition for Gibbs expectations and Connes’ classification of factors.  And today, Sergio’s student Roberto Longo has become one of the leading experts on modular theory.

I believe that my next encounter with Derek came in 1968, when he visited the ETH Seminar for Theoretical Physics during the summer, in the middle of my ten-week stay there as guest-professor. The high point of the afternoon at the friendly Hochstrasse 60, Zurich address involved discussion  around a table filled with cake and tea, organized by the institute secretary, Fraulein Rosemarie Hintermann. On these occasions Derek held forth with insight and his infectious humor, and as a result we renewed our friendship.

Derek and Marion moved to Bandol when Derek joined the group of Daniel Kastler. I was fortunate to have a long visit, living in an apartment in Cassis. I had interesting excursions to their home, and enjoyed the terrace, including one evening with Derek’s friend Tini Veltman.

Derek visited Harvard in the period 1971-2 as my guest, also as the guest of Shelly Glashow and Sidney Coleman. They three had taught together in a school organized by Feza Gursey in Istanbul, and Shelly was a friend with Derek at CERN. It was in Cambridge, Massachusetts that Derek collaborated on his famous “Lieb-Robinson bound.”  Their result is fascinating, as it can be regarded as an analog of the finite propagation speed that Glimm and I had just shown for the two-dimensional relativistic quantum fields that we had constructed.

After Derek left Marseille for Canberra, he made a tremendous effort to help build up mathematical physics in Australia. Derek hosted me twice in Canberra.  My first trip “down under” revolved around a January 1982 summer school. Elliott Lieb and I travelled there together, and afterward, Derek organized a grand tour around Australia for the two of us. During the hospitable time in Canberra, I got to visit Derek, Marion and their two girls frequently.

\begin{figure}[h!]
	\centering \includegraphics[width=0.25\textwidth]{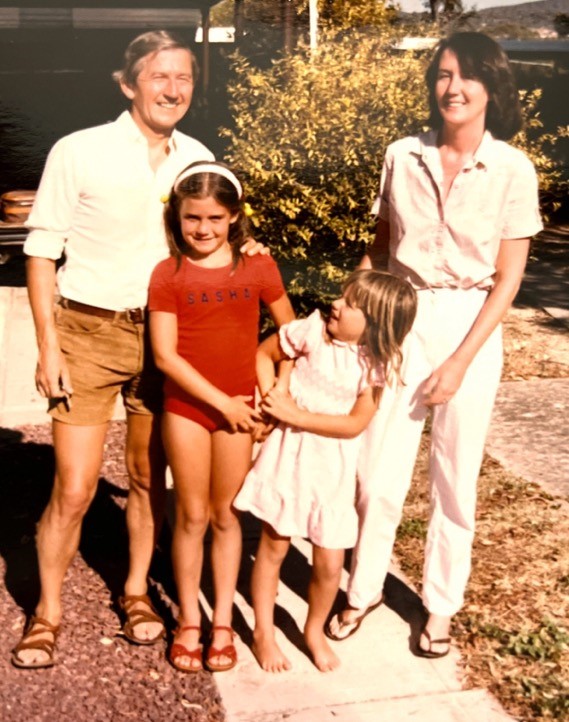}
	\caption{Derek and his family at his home in Canberra, January 1982}
\end{figure}

My second visit to Canberra in 1987 followed speaking at the annual meeting of the Australian Mathematical Society. The photos from these nice occasions chronicle the growth of Derek and Marion’s children.

\begin{figure}[h!]
	\centering \includegraphics[width=0.35\textwidth]{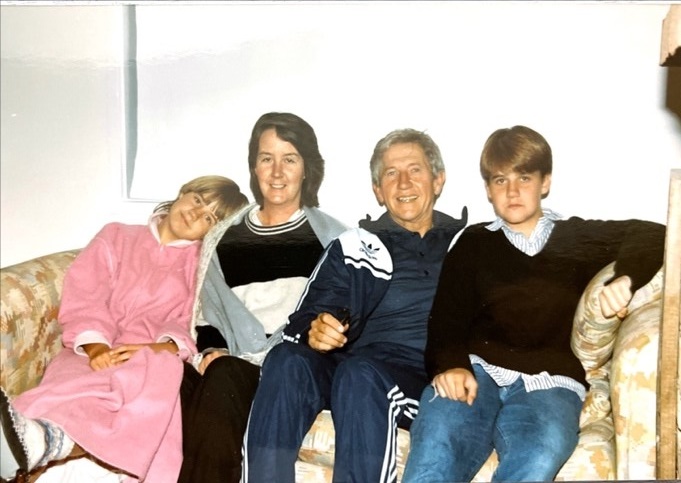}
	\caption{Derek and his family at his home in Canberra, 1987}
\end{figure}

Due to Derek’s life on the other side of the world, we did not interact as often as I would have liked. Our early interactions left a special imprint on me, including his humor and the unique timbre of his voice. What sticks in my mind are Derek’s stories of his unconventional friends, his love and penchant for cycling, his opinions about life in France, and the prospects for science in Australia. Sadly, now I only hear Derek in my mind; I miss his vibrant presence dearly.
%

\section*{Palle Jorgensen}

Over his life Derek led multiple collaboration teams. One research focus of collaboration inspired by Derek over multiple decades, starting in the eighties, was the theory of continuous one-parameter groups and semigroups on Banach spaces, leading in turn to new non-commutative $C^*$-algebraic structures.

\begin{figure}[h!]
 \centering \includegraphics[width=0.3\textwidth]{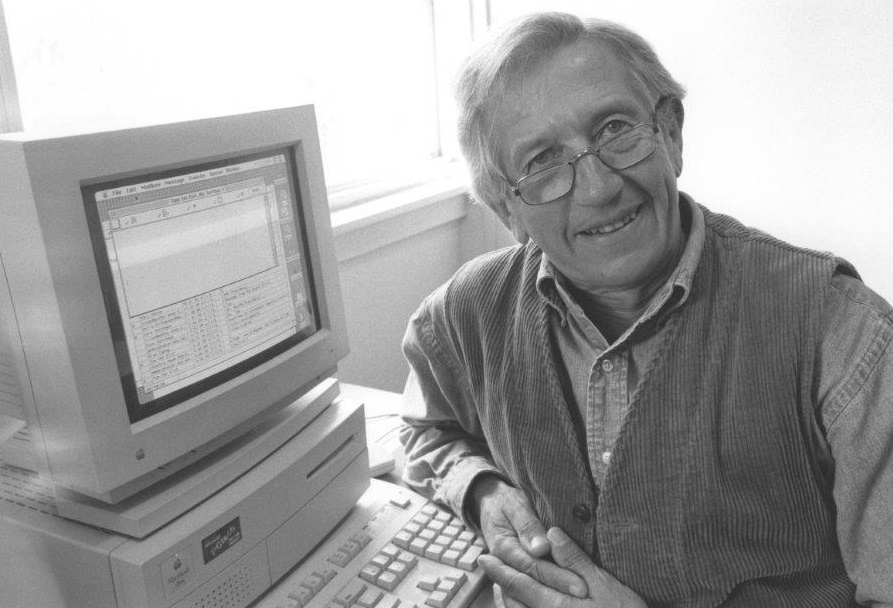}
   \caption{Derek Robinson at the computer, 1995}
\end{figure}

The collaborations started with annual research visits by Bratteli and me to Derek in Australia, during winters. A bonus with this arrangement for us visitors to Derek was that it allowed us to replace the cold winter in the North, with warm December-January Australia weather. Subsequently, the collaboration team expanded adding Charles Batty, Aki Kishimoto, Dai Evans, among others. Key themes  included an analysis of operator commutation relations, and their role in a unification of diverse areas in mathematics, in elementary particle physics, quantum-mechanical commutation relations, integration of Lie algebras, and unitary representations of noncompact Lie groups. A loss of a pioneer in mathematics leaves us with a void. In the case of Derek, for many of us, it was also a loss of a dear friend.
%
%

\section*{Aki Kishimoto}

While I was a Ph.D. student under H. Araki, Derek Robinson was the first specialist who gave encouragement via a letter  (undoubtedly responding to my mentor's request). This was cemented by my visit to Marseille  in October, 1977, a year later than originally planned.  This visit, while he was preparing to leave for Sydney the next year, lasted a year. Although my interest turned out to be confined to a narrow realm of operator algebras, mostly anchored to the book which he was then writing with Ola Bratteli \cite{BratRob} when I went there, his interest drifted away as time went by. But we managed to keep contact and occasionally meet until I received his last emails explaining his terminal health condition, adding that the paper he was revising then would be his last.

\begin{figure}[h!]
	\centering \includegraphics[width=0.4\textwidth]{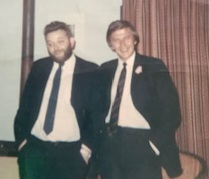}
	\caption{Derek with Best Man\\Tini Veltman at his wedding}
\end{figure}

When I visited him in Sydney and later in Canberra, his coauthor Ola was almost always present. Even after he found collaborators  in other fields, we talked about operator algebras related to mathematical physics in a broad sense. So he kept an interest in operator algebras from a physicist's point of view and shared his ideas with us. On my last visit to Canberra in 2006, we wrote a paper on flows. Although I met him a few more times after this collaboration, in Sapporo and Oslo, this was our last paper.  What I remember well during this period is his kind response to the messages (which must have been somewhat distraught) I sent to a few friends after the 2011 earthquake and Fukushima disaster in Japan. He offered his guest house as a refuge from the nuclear fallout.

I suppose there have been many mathematical improvements to the results in their book after it was published. But now they do not seem to be that important. In hindsight it was high time that Derek and Ola were finishing to write the book when I had a last chance to witness it and joined  some discussions on a couple of problems they came up with. Derek was quite proud of it as they could revise it as much as he wished.


\section*{Elliott Lieb}

It is half a lifetime since Derek and I worked together, or even saw each other, but I remember our discussions and friendship with clarity. Perhaps the last time we met was when my wife, Christiane, and I enjoyed a stay at Derek's and Marion's house near Marseille on our way to a summer school in Carg\`{e}se. In any case, I cannot forget our work in Cambridge Mass., leading up to our 1972 paper \cite{LR72} containing the bound on the speed of propagation of information in quantum spin systems.  The paper took over three decades to surface, thanks to Matt Hastings, who cited it in 2004. This is remarkable since it is about twice as long as it takes the cicadas we have just experienced in New Jersey to hatch.  It was fun doing it with him and I learned a lot in the process. I am grateful for that.

\begin{figure}[h!]
	\centering \includegraphics[width=0.4\textwidth]{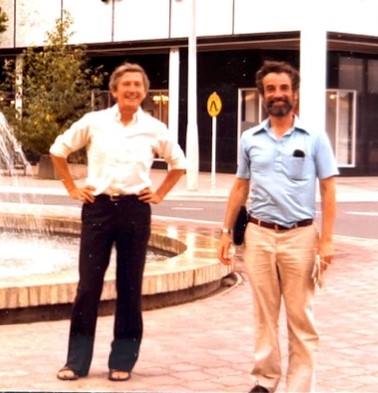}
	\caption{Derek with Elliott Lieb in Canberra, 1982}
\end{figure}

There are other ways in which the mathematical physics community and I are indebted to Derek. Among them there is the fantastic two-volume treatise with Ola Bratteli \cite{BratRob}, which established the foundations of our field. Equally important was Derek's foundational 1968 paper with Oscar Lanford on entropy \cite{LandRob}, and its conjecture of quantum Strong Subadditivity, which I learned about from David Ruelle. It was a bold conjecture because, as it turns out, many of the simple properties of classical entropy fail in the quantum domain, yet this very complicated property actually carries through from classical to quantum.

The conjecture kept me busy for several years, up to my 1973 paper \cite{LiebConv}, which contains the concavity of $A \to {\rm trace}\ e^{K + \log A}$. This, together with Beth Ruskai's
insight on the importance of the joint concavity of conditional entropy, $S(A,B)-S(A)$, formed the basis of our joint work \cite{LiebRusk}, which proved the fundamental Lanford-Robinson conjecture about entropy.

Those were glorious days in mathematical physics, generally, and Derek was one of its shining lights.


\section*{Bruno Nachtergaele}

To prepare for my official start as a graduate student in Leuven, Andr\'e Verbeure assigned me some `light summer reading' from the two-volume book by Bratteli and Robinson. I found the material fascinating and inspiring (even if a bit dense!). The book has been my constant companion throughout graduate school and until this day. It also was while a graduate student in Leuven that Mark Fannes pointed out to me the finite group velocity theorem by Lieb and Robinson. I talked a lot to Mark and he frequently alerted me to beautiful results in mathematical physics that seemingly randomly popped up in our discussions. So, I was well aware of Derek Robinson's work  but I cannot recall I ever met Derek in person. My only extended email exchange with him occurred in 2019 when he was reading the review paper I wrote with Sims and Young\cite{NSY}. He set us straight about some of the original history and pointed us to a follow-up paper of his we had been unaware of \cite{Rob76}. As my tribute to Derek's legacy, I will now discuss the subject of Lieb-Robinson bounds.

In non-relativistic quantum mechanics of interacting degrees of freedom, correlations between all of them are generated in any amount of time. At first glance, it appears that signals can be transmitted instantaneously or at least at arbitrarily high speed. In a $1972$ paper \cite{LR72}, Elliott Lieb and Derek Robinson showed that the dynamics of quantum spins on a lattice with finite-range interactions does, in fact, exhibit a bounded speed of propagation up to an exponentially small correction.

This result received some attention right away and it was recognized early on as a deep and important property of quantum many-body dynamics. It implied the existence of an approximate light cone which makes the non-relativistic dynamics look similar to relativistic quantum field theory of which the locality structure is a fundamental axiom. In spite of this, a series of important applications did not start to appear until more than three decades later. The citation history of Lieb and Robinson's paper, as shown by Google Scholar, hints at an event around $2004$ that ignited a steadily increasing interest. Up to that point, the paper was cited only a couple dozen times (so, on average, less than $1$ per year). Google Scholar now reports almost $1500$ citations, and that figure continues to grow at a steady rate of about $150$ per year.

Matthew Hastings published several papers in $2004$ that cite the Lieb-Robinson paper, most notably \cite{Hast1} and \cite{Hast2}. In the latter he gave a multi-dimensional version of the celebrated Lieb-Schultz-Mattis Theorem (LSM). In a eureka moment (which according to his recollection occurred while walking the Paris boulevards) he came to the realization that a finite propagation speed must hold for lattice systems with short range interactions and it was exactly what was required to complete his arguments. This was before becoming aware of the article by Lieb and Robinson with a proof of the property he needed. Another ingredient of Hastings' Lieb-Schultz-Mattis result is the Exponential Clustering Theorem for quantum lattice systems. Several proofs of such a theorem were given in the context of axiomatic quantum field theory decades before, but it remained unproven folklore in statistical mechanics. The proof of exponential clustering in theories with a massive vacuum state or in the presence of a spectral gap above the ground state uses locality in an essential way. In quantum field theory this locality stems from the finite speed of light (Lorentz invariance). In statistical mechanics, a proof had to wait for Hastings' discovery of how to use the finite propagation speed provided by Lieb-Robinson bounds as a replacement for the finite speed of light \cite{NS, Hast3}.

\begin{figure}[h!]
	\centering \includegraphics[width=0.30\textwidth]{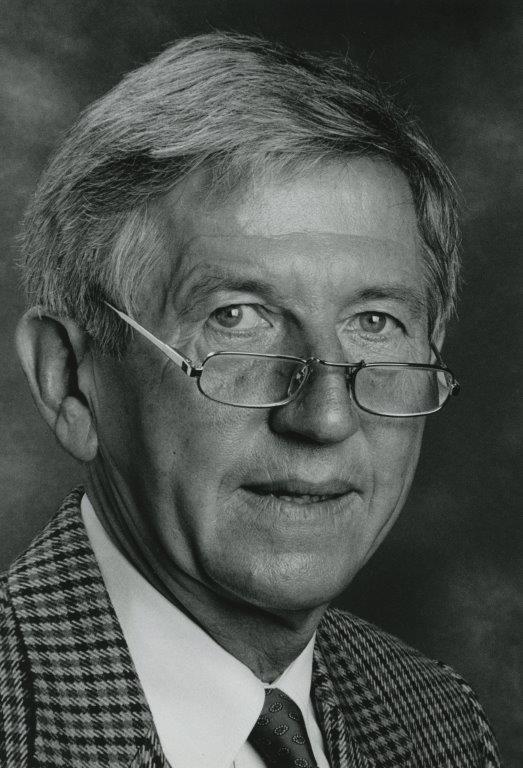}
	\caption{Derek Robinson, circa 1990}
\end{figure}

The Exponential Clustering Theorem and the multi-dimensional Lieb-Schultz-Mattis Theorem were the beginning of a sustained stream of applications. Lieb-Robinson bounds quickly developed into a fundamental component in quantum many-body theory and quantum information. Lieb-Robinson bounds are now established as a conceptual and practical tool in quantum theory. The vigorous pursuit of generalizations and improvements tailored to specific situations continues unabated.

%

\section*{Heide Narnhofer}

Presumably Derek Robinson did not know about his influence on my scientific life, but it is a fact and I am happy that now I can express my gratitude.

I first met him when I was a student and he gave a talk in Vienna, how in the language of "Local quantum theory" the thermodynamic quantities energy, entropy and pressure can be obtained. This theory was fairly new at that time and Derek was one of its pioneers. The demonstrated clarity how physical quantities can be expressed as mathematical objects and on this basis how rigorous considerations can lead to new expressions and realizable facts impressed me very much and I felt immediately that this is the language in which I could approach physical problems.

\begin{figure}[h!]
	\centering \includegraphics[width=0.3\textwidth]{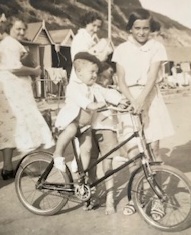}
	\caption{Derek as a young cyclist in the 1930s}
\end{figure}

Three years later I was a postdoc at IHES in Bures and participated in a workshop of several researchers in mathematical physics, one of them Derek. No talks were given, instead someone got up, starting with: "Recently I have thought about this problem.." and then the discussions and suggestions started. Derek spoke about the influence of boundary conditions, offered partly results and also the lack of controlling the independence on Dirichlet and Neumann boundary conditions in the thermodynamic limit. R.B. Griffiths reacted, and I could follow, how out of intuition, confidence in intuition, some experience in mathematical modeling and the exchange of ideas new ideas evolve and results emerge. In my following research I always tried to copy this combination of strong belief mixed with criticism without being discouraged by it.

My last experience with Derek was in Marseille, when I asked him how to understand why pure thermodynamical states have more physical relevance than their mixtures. Derek offered immediately the explanation given by their stability under the perturbation of their dynamics. There only remained the task to polish the argument. Years later I gave another explanation based on the idea of collapse in the framework of quantum history.

This corresponds to Derek's favour: he turned more and more to mathematics and especially to the beauty of algebraic theory and its connection to differential theory, and he searched in this area for results. The theory of collapse that is based on randomness where probability theory is needed to find some structure was less attractive for him. I am sorry that we never had the possibility to discuss which approach reflects better the reality. This would have been very fruitful, because Derek was very open minded, and it was easy to come into contact with him. He never gave the impression that one might ask stupid questions but instead was encouraging and stimulating.


\section*{David Ruelle}

{\em A conference was organized in Derek’s honor shortly before his death.  David Ruelle sent the following letter to Derek at the time which we include with Ruelle’s permission.}

As you pointed out to me, our joint paper is still quoted after all these years. Of course science goes on, but our work took place at a crucial time and has left its imprint.

We have met over the years at various places, notably Bandol with Daniel Kastler when you seemed permanently established in France, and later Canberra where you invited me to visit the Australian National University. But I think that the first time we met was at the ETH in Zurich. This is when, after the death of Wolfgang Pauli, Res Jost congregated a little group of visitors at the ``Theoretische Physik'' first at Gloriastasse 35 then at Hochstrasse 60.

Let me remind you particularly of your Zurich period. Apart from you and me, young people who were around were Othmar Steinmann, Gianfausto Dell'Antonio, Klaus Hepp, Peter Curtius, Huzihiro Araki, Walter Hunziker, and others. We went for lunch at the Vogelsang where my wife Janine joined us. We could drink Apfelsaft there or alkoholfreies Rum Punch, and discuss forever on the future of physics. I also remember a party in an anti-atomic shelter of the Swiss army by the lake of Zurich. Those were the good days...


\section*{Mary Beth Ruskai}

In the 1960's, mathematical physicists seeking a firm foundation for statistical mechanics found that the framework of operator algebras was well-suited to this task.  When I was a young postdoc in Geneva in 1970-71 two of the leaders in this area, David Ruelle and Derek Robinson, gave Troisi\`{e}me Cycle courses in Lausanne with Derek's lectures focussed on the thermodynamic pressure.

Earlier, Derek co-authored a pair of important papers on entropy.   The first, with David Ruelle, studied classical systems, while the second, with Oscar Lanford, studied quantum systems. David suggested that I extend the results in the latter to von Neumann algebras with a well-behaved trace.  Derek, then at the University Marseille, invited me to come for a week to coincide with a visit by Oscar.  He arranged for me to stay at a small hotel near his home in the lovely town of Bandol on the Mediterranean.  I fondly remember  meeting each afternoon at one of the cafes along the beach to discuss mathematics.

\begin{figure}[h!]
	\centering \includegraphics[width=0.4\textwidth]{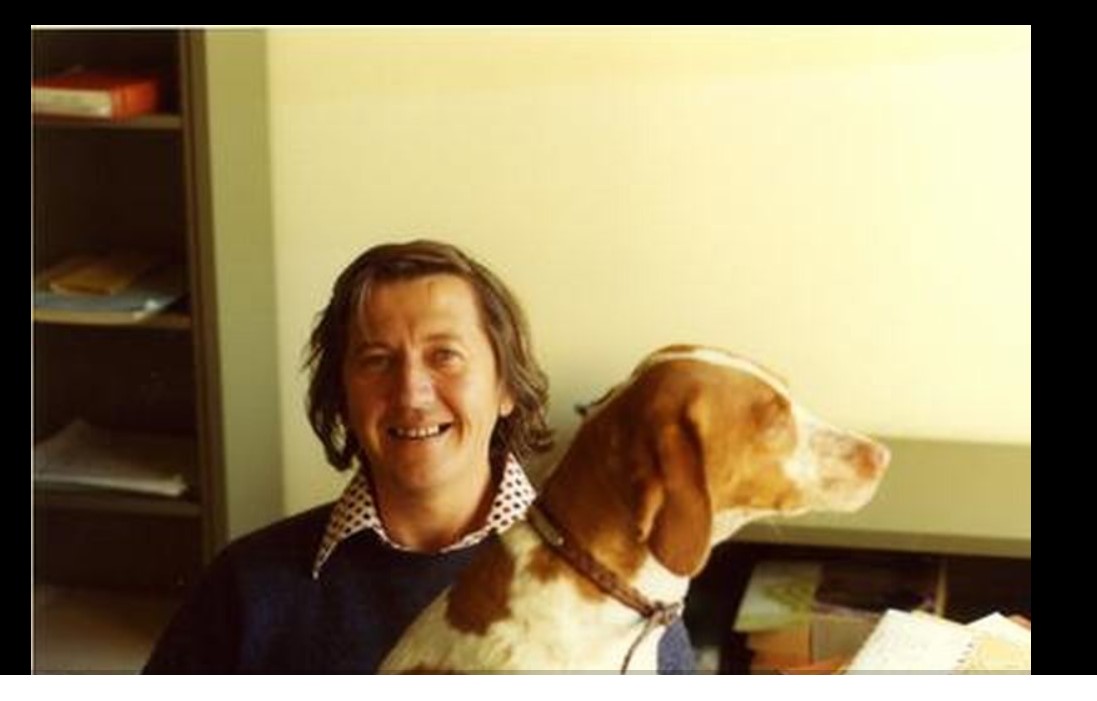}
	\caption{Derek at Berkeley about 1975}
\end{figure}

While doing this work, I realized that if  the conditional entropy $S(\rho_{12}) - S(\rho_1)$ were concave it would imply a conjecture at the end of the paper by Lanford and Robinson.   It may seem strange to suggest that the difference of two concave functions would be either concave or convex.     However, Elliott Lieb proved that the conditional entropy is, indeed, concave and this was an important ingredient in my work with him proving the strong subadditivity of quantum entropy.

After Derek left Marseille we didn't cross paths again until 1997  when we met at the International Congress of Mathematical Physics  in Brisbane.  He was  pleased to hear that the work on quantum entropy he began with Lanford in 1968 was now having an impact in quantum information theory.


\section*{Adam Sikora}

The first time I met Derek was long ago, in October 1995. We met in Canberra, on the Australian National University campus, in a building that no longer exists. I came to Canberra on my first postdoc position, funded by Derek's ARC grant. It was easy to start working with him. Derek's research expertise and interest in mathematics were always extensive. I did not have to learn that much of anything new - we just started to talk to find a topic appealing to both of us. And it was always a pleasure to meet and chat with him. He had an excellent sense of humour and lots of exciting things to share.

In our first joint work, we studied Riesz Transform in the setting of the Lie group. And it was the beginning of a quarter of a century of collaboration and friendship. After Lie groups, we investigated operators with periodic coefficients. Next, we were interested in degenerate elliptic operators. There are still many topics that I would like to discuss with Derek today, but sadly this is not possible anymore.

Derek always wanted to learn new things. Even in his sixties or seventies, he was not afraid of starting a study in a new research area and learning new things. The evolution of his research expertise is impressive and broad: from physics to C*algebras, Lie groups, Riesz transform, operators with periodic coefficients, degenerate elliptic operators and many more.

When it was clear that Derek was terminally ill and his prognosis was not good, I had a chance to talk with him a few times over the phone. Derek told me that he had submitted a new paper to Journal of Functional Analysis, that he was frail and struggled to revise the manuscripts, and how happy he was when it was accepted. It made me realize how significant mathematics was for him. He was very tired and not able to chat for too long, but we were able to discuss many other things. I am surprised by how important it was for me to talk with him and to thank him for our collaboration and all the other things he did for me.

Derek had the ability and wisdom to cherish every day of his life – how to use any moment and enjoy it. As we all know, he loved cycling, but when cycling became too risky for his age, he took up race-walking to appreciate good weather and stay fit and healthy. Derek loved to be active, which was always a part of his personality. First of all, he always wanted to do some new mathematics. Make some exciting calculations, write and explain them to others and discuss them over coffee with friends. Derek used and enjoyed every day of his life as much as possible.


\section*{Barry Simon}

I first met Derek in early in 1973.  Freshly tenured at Princeton, I spent AY 1972-73 visiting three great centers of European mathematical physics: the fall at IHES near Paris, the spring at ETH in Zurich and the winter at what was jokingly called the free university of Bandol.  This was the name given to the remarkable group that Daniel Kastler had put together at the Universities of Aix-Marseille and Toulon and mainly at the CNRS in Marseille where I was officially visiting.  Two of the other senior members of the group were Derek and Alex Grossman.  Daniel and Derek lived in the village of Bandol about 50 km from CNRS and Alex lived in the quaint fishing village of Cassis half way between which is where my wife and I stayed.  Much research was done informally in the villages, hence the name.

Derek was lively and outspoken.  I was somewhat shocked by his references to {\em damned frogs} although I eventually learned about his travails  with French bureaucracy over arranging foreign visitors and building his house and understood this usage better.  In any event, Derek did it with such charm that it seemed all right in any event.

My second and more intense interaction with Derek was when he invited me to visit in the summer of 1983 or perhaps I should say the winter since the visit was to Canberra where he’d relocated.  The visit almost didn’t happen.  My youngest son was born in Dec. 1982 and by the time we were able to get a birth certificate and use it to get a passport for him, it was the end of April before we sent the visa application to the Australian consulate.  The visa officer was making unreasonable demands about a medical form so, in the era before email and with expensive international phone calls, I sent a telex to Derek explaining that I might not be able to come.  He called and told me to hang tight.  Two days later, the formerly officious visa officer called: ``Sir, I am anxious to issue your visa but I need you to return your passports.’’  ``What about the medical form?’’ ``Oh you don’t need that, sir’’.  Derek had performed a miracle.

\begin{figure}[h!]
	\centering \includegraphics[width=0.4\textwidth]{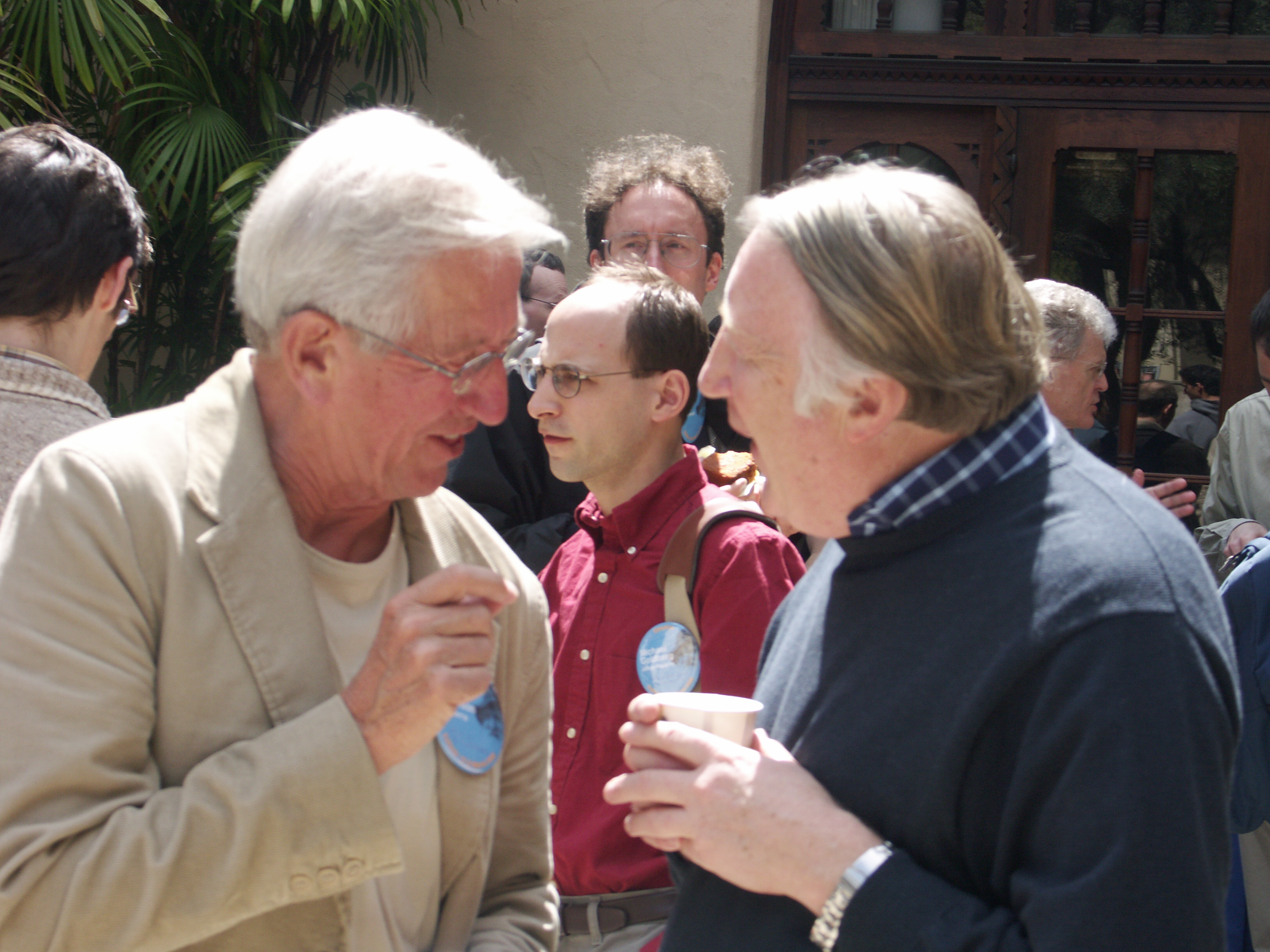}
	\caption{Derek with Aubrey Truman\\at Simonfest, 2006}
\end{figure}

It was a very important visit for me scientifically.  Michael Berry was visiting the physics department and told me about some recent work he’d done which led to my paper on what I called {\em Berry’s phase} surely my most famous work in the physics community.

Brian Davies was also visiting and we worked hard on some results involving a new concept that needed a name.  We found a notion which was stronger than the notions of \emph{hypercontractive semigroups} and \emph{supercontractive semigroups}. I remember the meeting we had in Derek’s office trying to figure out a term stronger than ``hyper’’ and ``super’’.  After we joked about ``super-duper’’, Derek suggested ``ultra’’ and so was born \emph{ultracontractive semigroups}, a terms with almost 10,000 hits on Google!

While our research interests weren’t distant, they weren’t really close except for two lovely things that Derek had in his book on the Thermodynamic Pressure~\cite{RobPress} which involved my focal interest of Schr\"{o}dinger operators.  In our work on Thomas Fermi theory, Elliott Lieb and I used a variant of a technique of Weyl which Andr\'{e} Martin had used to get the large coupling limit of the number of bound states of a Schr\"{o}dinger operator but we learned Derek had it independently.  And Derek had proven a conjecture of Kato on monotone limits of forms earlier than I had.

In recent years, Derek and I frequently swapped emails – exchanges that were illuminating and often fun.  I miss him.


\section*{Tom ter Elst}

Almost unnoticed Derek Robinson fulfilled major administrative responsibilities.  For example, he was the Chairman of the Board of the Institute of Advanced Studies at the ANU from 1988 till 1992. In order to stay research active, he decided in 1988 to write the book \cite{RobLie}. On my arrival day as a postdoc in Canberra mid 1990 he gave me the first two chapters of the manuscript with the message that it was the final draft and that I certainly would find a topic to work on. I was deeply impressed how well it was written. It contained a wealth of information and many things to work on. The `final draft' still got two revisions.

\begin{figure}[h!]
	\centering \includegraphics[width=0.4\textwidth]{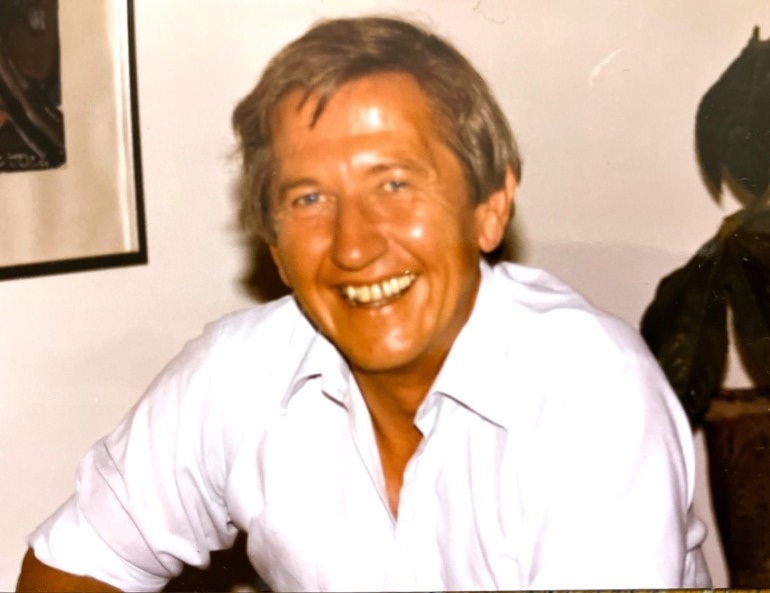}
	\caption{Derek in Canberra, 1982}
\end{figure}

The book covered (possibly higher-order and complex) strongly elliptic operators and sublaplacians on Lie groups. Derek used a Nash inequality in order to obtain Gaussian bounds for the kernel of the semigroup associated with elliptic operators, extending Langland's PhD results. Together we managed to extend many theorems to the setting of possibly higher-order complex subelliptic operators on a Lie group, even allowing for different weights in different directions, as happens for positive Rockland operators on graded Lie groups. We next describe a striking example. Let $L$ be the left regular representation of $SO(3)$ in $L_2(SO(3))$. Further let $A_k$ be the infinitesimal generator of the natural one-parameter groups. Then the ordinary Laplacian is the operator $A_1^2 + A_2^2 + A_3^2$ and a sublaplacian is $A_1^2 + A_2^2$. Derek investigated these operators. Now take as example $H = A_1^4 - A_2^6$ or $H = A_1^4 - A_2^6 + A_1^2 \, A_2^3$. Later, Derek and I proved that operators like $H$ are also the (minus) generator of a holomorphic $C_0$-semigroup which has a kernel satisfying Gaussian type bounds, including all its higher-order derivatives. This even provided a characterisation of these weighted subcoercive operators.

A different challenge is to obtain asymptotics for the semigroup and its kernel for large time. For operators with periodic coefficients a major step was obtained with techniques of homogenization theory.  The details were written in another book \cite{DER03}.
After that Derek got more and more interested in degenerate elliptic operators and the relation with Hardy inequalities \cite{Rob21A}.


Michael F. Barnsley is Emeritus Professor at the Mathematical Sciences Institute, Australian National University.  His email address is michael.barnsley@anu.edu.au

Alain Connes is professor Emeritus at College de France on the chair ``Analyse et G\'{e}om\'{e}trie'' and at IHES on the ``L\'{e}on Motchane chair''. His email address is alain@connes.org

David Evans is an Honorary Distinguished Professor at Cardiff University. His email address is evansde@cardiff.ac.uk.

Giovanni Gallavotti is Emeritus at Universit\`{a} ``La Sapienza'', Roma. His e-mail is giovanni.gallavotti@roma1.infn.it.

Sheldon Lee Glashow is Higgins Professor of Physics emeritus at Harvard University.  His email address is slg@bu.edu.

Arthur Jaffe is the Landon T. Clay Professor of Mathematics and Theoretical Science at Harvard University and a past president of the American Mathematical Society.  His email address is:  arthur\textunderscore jaffe@harvard.edu

Palle Jorgensen is a Professor of Mathematics at the University of Iowa, USA. His email address is palle-jorgensen@uiowa.edu.

Aki Kishimoto is Emeritus Professor, Department of Mathematics, Hokkaido University, Sapporo, Japan.  His email address is akiksmt@r3.ucom.ne.jp

Elliott Lieb was a professor at The University of Sierra Leone, Yeshiva University, Northeastern University and M.I.T. and currently professor emeritus at Princeton University.  His email address is lieb@math.princeton.edu

Bruno Nachtergaele is a Distinguished Professor of Mathematics at University of California, Davis. His email is bxn@math.ucdavis.edu.

Heide Narnhofer is a retired professor of mathematical physics, Faculty for Physics, University of Vienna. Her email is heide.narnhofer@univie.ac.at.

David Ruelle is Emeritus Professor at Institut des Hautes \'{E}tudes Scientific (IHES).  His email address is ruelle@ihes.fr.

Mary Beth Ruskai is retired.  Her email is address mbruskai@gmail.com.

Adam Sikora is an Associate Professor of Mathematics at Macarise University. His email address is adam.sikora@mq.edu.au.

Barry Simon is the IBM Professor of Mathematics and Theoretical Physics, Emeritus at Caltech.  His email address is bsimon@caltech.edu.

A.F.M. (Tom) ter Elst is a professor of mathematics at the University of Auckland. His email address is terelst@math.auckland.ac.nz.

\textbf{Photo Acknowledgements}

``Derek Robinson (1935-2021)'': Photo courtesy Michael Hood from Faces of Science Exhibition \\

``Derek and Marion in Bandol'', ``The House in Bandol'', ``Derek..Baton Rouge Conference'', ``Derek and his family$\dots$1982'', ``Derek and his family$\dots$1987'', ``Derek with Elliott Lieb$\dots$'', ``Derek in Canberra, 1982'': Photos courtesy Arthur Jaffe \\

``Derek with his cycling medal, 2002'', ``Derek on His Bike'', ``Derek with Ola Bratteli'', ``Derek as a speed walker'', ``Derek with Best Man Tini Veltman at his wedding'', ``Derek as a young cyclist in the 1930s'': Photos courtesy Marion Robinson \\

''Lysebu,$\dots$'': Photo courtesy Eric Christopher Bedos

``Derek with His Bike', 1995'', ``Derek Robinson at the computer, 1995'', ``Derek Robinson, circa 1990'': Photos courtesy ANU Archives \\

``Derek at Berkeley about 1975'': Photo courtesy George Bergman \\

``Derek with Aubrey Truman at Simonfest, 2006'': Photo courtesy Barry Simon

\section*{Derek's Autobiographical Note}\label{D1}

My father was born in rural Ireland quite close to the border in County Monaghan. He left school at the age or 10 or eleven for farm work and opportunities for him were very limited. As he grew older he looked for better options and the one he chose was to join the British Army. In fact he joined it twice, the first time before he was legally old enough but this was discovered and he was dispatched home. He joined again, this time successfully, about 6 months later.

He spent 3-5 years with the army in England and then left the army but stayed as a permanent member of in the army reserve. This was important to our family as when the 2nd World War was looming, before its official beginning, he was called up to rejoin the army and the whole family was relocated to the Larkhill army base on Salisbury Plain. We spent the war years in the Larkhill camp and this was where I had my early education – in a tin hut on Salisbury Plain – with a mixture of teachers, all female, mainly retired women teachers or wives of serving soldiers who were helping out. Initially all the children were sorted by age into classes, but it was found that with the limitations on teachers and space this was not ideal so they sorted them by ability and I ended up in the top class with two years to go before the official end of primary school!

My mother left school at the age of 14 and my memory says she was working in retail and living in Croydon. I am not sure how they met.

I have an older brother who is nearly two years older than me. When we got to Larkhill I would have been 4 and he would have been nearly 6. When the war ended, because we had been the first to be called up, we were among the first to get out. My father was demobbed in '45 and given the traditional free suit and we moved from Larkhill to Bournemouth were we had lived for a few years before the war. I still had a year of primary education to go before secondary school and my father was very keen I should get the best education for that year and I remember him going around to all the primary schools within cycling distance to try and check out which schools were good and which would take me. One school annoyed him greatly because he told the teacher he had this son who had done very well in this army school and this was dismissed with – oh, there are so many bright guys in these rural schools but now you are in the big city! But of course I carried on doing what I did and by the time I finished the final primary year I was 5th in the class and ready to move on to the grammar school. I must say that during this whole process of education I never really understood what was going on with the grammar schools and the other schools. I had this vague idea that not going to the grammar would result in my becoming a carpenter.

Things just progressed from there. I was normally near the top of the top stream in my class. It was built on a natural affinity for mathematics and if you do well in mathematics you are half way there. I could do well in other subjects but I was too damn lazy so I often didn't try in them.

At the age of 14 or 15 I developed what is now called mononucleosis, it was unknown and unnamed and untreated and I ended up spending 6 or 7 months in bed. About that time my brother was choosing his career. He had been fascinated by boats from the age of 2 and was determined to go to sea in some capacity. My father found that various standard commercial companies offered cadetships and my brother went with British India - the East India Company – where he went on to have a long and illustrious career. This influenced me as its natural to think about doing the same as your older brother so I considered a seagoing career and I looked to see what was involved in getting accepted. Of the required academic qualifications I found that mathematics was important and got you a step up in the process so I looked into what mathematics they wanted and found it was this mysterious thing 'calculus'. So I decided to teach myself calculus and so for about 6 months I worked away and learned calculus. As far as joining the navy was concerned this was wasted effort as I gave up on that idea but as far as school was concerned this was a bit undermining for the school! I knew all sorts of things that others didn't. By the age of 16 or 17 when the school was teaching calculus I was way ahead or certainly sufficiently ahead to make a difference. So this was an example of my older brother's career influencing mine.

It gradually became clear to me that the aspiration of the grammar school was to put as many people into university as possible and this was fine by me. The only other thing I was good at was cycling and given the choice of professional cyclist or mathematics at Oxford I chose mathematics. Curiously enough the 2021 women's Olympic cycling road race champion is a professional mathematician and an amateur cyclist. I had no idea what a university actually was, I really seriously thought there were only four universities : Oxford, Cambridge, Harvard and Yale and that was it. But I realised that getting to university was the next stage and since my father always wanted me to keep on to the next stage, I kept going and my progress though sixth form at grammar school was fairly straightforward. I don't remember the details but based on advanced level exams I got offers from London and Bristol. It became apparent that if you wanted to get to Oxford or Cambridge from a 'normal grammar school' more was needed than Bournemouth Grammar had provided at that stage. I realised later when I reached Oxford that the people from Manchester Grammar School had had quite a different level of preparation! I decided to stay on at school for an extra year to do the entrance exams for Oxford. This ended up being half a year as I started working to earn some money. Luckily there were three other guys also applying to Oxford, two in chemistry and one other in mathematics and we all ended up getting offered places. Entrance to Oxford then was via the colleges and I remember staying in Wadham College for a period of two days when there must have been some sort of written exam and we were extensively interviewed with the outcome I was offered a place which I accepted – again not knowing what I was letting myself in for and knowing very little about the college. There was a small hiccup that I needed to pass the ordinary level in Latin which I had failed first time around. This I succeeded in passing on the second try even though they had changed the 'set book' between my first and second attempts. All I can remember from this is that all of Gaul was divided into 3 parts – which struck me as sensible way of doing things!

What I found out when I got to Oxford in 1954 was that the head tutor at Wadham College who influenced everything for the undergraduates was an applied mathematician and a pretty classical applied mathematician so this influenced my undergraduate work towards applied mathematics. But the one thing I really appreciated was there was no obligation to do anything and I spent a lot of time doing nothing! I started going to first year lectures down the main street but I realised I wasn't learning very much that way and the one hour in applied mathematics with my tutor was getting me further ahead so I thought why bother getting up early and going to lectures when I could do more by myself. So I got into the habit, very early on, of working by myself. I also got into the habit of not working but rowing and drinking beer. This had an influence on my career because Oxford only measured our progress with exams at the end of first year, called mods, and then finals at the end of the third year. The mods were important as the top students in mods got honours, you could make your name by doing well in mods. I went into the mods exam expected to achieve first-class honours but this didn't happen as I fell asleep during the exam due to my various social activities! So that left me a poor old commoner until finals at the end of the third year.

The third year you had special subjects and I was naturally pushed into applied maths subjects. There were two subjects that appealed to me: one was hydrodynamics and I don't know why that appealed as I quickly lost interest in it and the other was quantum mechanics. That was where I really started getting interested in something really interesting. I realised that at end of the third year there were two options, either to get a good enough degree I was offered a job which was exempt from conscription or to spend two years in the army. I didn't like the second option so I worked quite hard on quantum mechanics and in a surprise to most other people, at the end of the third year I got a first!

That left the situation sort of hanging because the only person who had thought I might do well was myself so I had not been able to apply for a doctorate positions without support from others. But we had a new head of department in theoretical physics who I had run into due to doing quantum mechanics, Willis E Lamb Junior of the Lamb shift, and he was setting up a one year diploma course which was something I could apply for so I had applied for that. Now once you have a good first there are all sorts of people from all sorts of places saying come work with us, we have a scholarship. Most of those people were theoretical physicists. I got tempted and lured by various different places but I was also offered a scholarship from Oxford because I had done better than they expected so I went into doing physics at Oxford - that was how the cookie crumbled! It was a blessing and a curse as the only openings they had in theoretical physics were in nuclear physics, down and dirty stuff, so I ended up having to work though this stuff in which I had absolutely no interest but that led me to spend time reading around and so I became interested in quantum field theory and particle physics. Indirectly that influenced my postdoctoral situation when I went to Zurich where I ran into many more people interested in these topics.

Within the nuclear physics people there was a guy Denys Wilkinson, a bright leading light who was probably only 10 years older then me and he had just been appointed professor of experimental cum theoretical nuclear physics cum particle physics and he ran informal evening talks given by internal people or people from nearby. That was where I first met Abdus Salaam as he came up from London and gave a talk. I went along to these events unwittingly and unknowingly that Denys was going to suddenly spring two papers on me: one on particles and one on anti-particles, there was this anomaly of their difference in behaviour. I said this were interesting and he said yes, how about giving a talk on it in three weeks time!! So I suddenly found myself doing something I was partially interested in and I think I gave two or three talks. It had a longer term influence on my career as he was head of the committee which gave me a postdoc position in Zurich some years later! So I think it was self interest both on his part and on my part.

Another person who worked down the road in Harwell was John Bell of Bell's theorem etc. I guess we knew each other in the sense we talked in groups after seminars and knew each others names but he was about 5 years older than me and we never got close. He turned up later in CERN and was a close friend of a group of other people I knew.

So I guess the next question is how I turned from a down-to-earth nuclear physicist into a mathematical physicist.

During my doctorate, the head of department, Willis E Lamb Junior who had won the Nobel prize about three years earlier was trying to bring us up to modern standards. One thing he did was have an allocation of money so all the graduate students could go to at least one international meeting. I was very lucky to convince them I needed to go to Naples for a summer school in field theory and particle physics which interested me a great deal. There were number of people there who were influential on my later career. One lecturer who impressed me was Res Jost, he was a professor in Zurich who had replaced Pauli who had died prematurely about a year earlier. Later on I ended up going to work for him in Zurich on a NATO fellowship. It was odd being in Zurich and working for NATO but it meant I could work there and cost them nothing but it caused questions and puzzlement as Switzerland is definitely not a NATO country! Another important contact I made in Naples was Tini Veltman who ended up getting the Nobel prize in 1999 for renormalization theory. Someone else who came out of the woodwork there was David Ruelle. David Ruelle is actually Belgian and the year I was in Naples he was cleaning up his PhD at Brussels but working in Zurich. Later when I got to Zurich he was one of the bright guys there and he and Michael Fisher were two of the leading lights in deducing exact results in statistical mechanics. So this meeting in Naples was important for me in lots of ways. I don't remember the names of all the other people in Naples but there was a large Dutch group around Tini.

Naples was about half way through my PhD and there was a similar summer school in 1960 immediately after the end of my PhD which was held in Edinburgh where I again met Tini. So we redid what we had done in Naples which was to form a wine club. We drank wine and spent our evenings trying to piece together the lectures. During the Naples event we had really not understood very much at all but by the time a year had gone by and we got to Edinburgh we were a bit more astute. We actually had quite some fun as we could bate some the lecturers in Edinburgh on some of the stuff we could see more clearly than they could! We refined the basis of the club down to 4 people: it was Tini and myself, Nicola Cabibbo of the Cabibbo angle and Shelly Glashow of the standard model, so 2.5 Nobel prize winners! I kept in contact with Tini all his life and I visited him often and he was our best man when we got married in Marseilles, I also kept in touch with Shelly. Now of course what they ended up doing and what I ended up doing was quite different but I understood more or less what they were doing and for some years I followed what they were doing without them following what I was doing which gave me a broader perspective on things and shows the importance of going to these summer schools in Naples and Edinburgh.

Tini had done his service with the army in Holland and he did his PhD with L\'{e}on Van Hove. Tini was a rather aggressive young man and he gave L\'{e}on Van Hove a very hard time and they were always arguing. It was a very strange situation for a relationship with a thesis advisor but it all got sorted out in the end. One interesting thing was that Tini ended up sharing the Nobel prize with a student of his, Gerard 't Hooft. Now Gerard 't Hooft's was not a new boy in physics because his uncle was van Kampen who was a leading figure in solid state, many-body theory. The idea was that 't Hooft would do his PhD in Utrecht with van Kampen but Tini undermined the whole thing by convincing 't Hooft he should work with him! So all sorts of tensions developed and despite being very productive this ended up with 't Hooft and Tini having a very 'distant relationship'!

I spent two years in Zurich (academic year 1960-1962) where I was part of a group of about 5 postdocs of which David Ruelle was the outstanding one. In the second year Araki came for 3 months and he gave lectures on Von Neumann algebras which I knew nothing about at that point, in fact none of us knew anything about them , but we were all impressed they had something had to do with field theory. Araki had come from Princeton where he was working as a postdoc with Rudolf Haag who was really the man behind a lot of these things.

When I had been in Edinburgh at the summer school I met a Canadian, Dave Jackson who was a colleague of Rudolf Haag at Illinois. Dave was not really interested in field theory but he was favourably impressed by me, Tini and our little group and when, later on, I applied for a job at the University of Illinois at Urbana-Champaign Rudolf received a reference directly from Dave. Rudolf hired me immediately which was both a good thing and a bad thing. It was a bad thing because it was one of the my hardest times. Araki and a group of others at Illinois were all way ahead of me and it took me a long time to catch up, but I eventually did. I was there for the academic years 1962-64.

After Illinois I went to Munich for eighteen months as that was where Haag was on his sabbatical leave from Urbana-Champaign. We had a couple of other people around who were interested in similar things and we were having a general talk one morning when there was a knock on the door and it was Sergio Doplicher and his cousin Guido Fano from Rome who had heard Rudolf was in Munich and had come to visit him on the hope he would talk with them. Now Sergio was already quite well advanced in these things as he was a smart cookie so Rudolf was impressed by him and he joined the small, adhoc group in Munich.

At the end of 1965 Haag went back to Illinois. This was the middle of the academic year so I spent 3 months in Paris at IHES with David Ruelle and 6 months in Marseilles with Kastler. Kastler also invited Doplicher to Marseilles so we spent some time together and the three of us wrote a joint paper which was maybe the first time I really started to believe there was light at the end of the tunnel, that I was going to get somewhere. During this period I had the advantage that I had many other areas of interest: quantum spin systems, Ising models etc. So while they were doing pure C* algebraic things and maybe hoping to find a model that modeled them I was also looking at the models, but not worrying about whether they immediately satisfied all the properties, leaving that as a problem to be solved – so I was working in my own little broader scene.

After Marseilles and Paris I went to CERN in Geneva for the academic years 1966-68. I'd applied to go there two years earlier but the only guy who could have been interested in me completely overlooked my application! When I applied again a year later he apologised for not having even noticed the earlier application. It was during this period I met Marion. First I met a NZ diplomat, Hugo Judd and he knew Marion's father who was working at the International Labour Organization and that led to a group of us going over to Marion's family house for a party. Marion had left NZ when she was only three as her father was appointed to a job in Canada and then Geneva. So she had been raised as a diplomatic child in Geneva, attending the international school and traveling back to NZ every two years. Not flying at that date of course, the journey by boat taking 6 weeks each way which lead to her disliking travel immensely!

In CERN my initial interests were in particle physics and I gave a series of lectures at a summer school in Istanbul. Shelly Glashow and Sydney Coleman were also lecturers. My lecture notes were distributed by CERN but my copy has gone missing with time. Despite this appropriate start I then spent much of the time in Geneva examining equilibrium properties of classical and quantum mechanical systems. For example, the uniqueness of equilibrium states at high temperatures and nonuniqueness at low temperatures. My main interest was in quantum systems, that is, generalizations of the Heisenberg model. The latter is, of course, a direct quantum equivalent of the classical Ising model. I managed to prove that the dynamics of the infinitely extended model and a broad class of generalizations was given by a strongly continuous one-parameter group of automorphisms of the algebra of spins. This was, in fact, not difficult. The only real problem was to overcome the general scepticism that the result was incorrect or only true in very special circumstances such as the Heisenberg nearest neighbour model. My good friend David Ruelle was a major sceptic who expressed both doubt and astonishment when I told him of my general conclusion.
It must be emphasised that during that period physicists were generally sceptical about the use of properties of infinite systems to understand the properties of very large systems. The use of asymptotic approximations was widely used without question in many areas but strongly resisted in statistical mechanics. I was often beset by irrational tirades against my introduction of infinite systems.

1968, after two years at CERN, my post-doc period was over. I had several offers of tenured jobs, in physics and mathematics and I was offered posts in Austria, Holland California and France. One of the results of the student uprising in 1968 was to make it possible for non-french people to have permanent university appointments, in fact permanent public service jobs, so Marion and I decided to go back to Marseilles where we stayed for 10 years. It was a very busy time!

The university was expanding and a new faculty was being created to the east of the city. I was appointed to this new faculty of science which still had to be created and this required a lot of work structuring the departments, the institutes and the Council. For my sins I was appointed to the transitory council which inherited many of the features of the revolutionary bodies of 1968, long, argumentative meetings were the order of the day.

One surprising attribute of the new faculty was the students(!) insisted the lectures should be in block form. Each lecture topic was covered by lectures on four days a week for three months of the year. Thus three full topics were covered each year. This suited me no end as I was essentially free to do research for 9 months each year.

Trond Digernes had just got his PhD in LA and he came to do a postdoc in Marseilles under a two-year Norwegian scholarship. One morning I went in and there was very quiet guy there – Ola Bratteli. I'd asked Trond about a problem on unbounded derivations and he said he'd asked Ola and Ola thought I was wrong which excited me greatly. Much more exciting to be wrong than to be right! Ola started explaining but he really was very quiet but the result was we spent half the day talking and this led to us working together for a while.

I had in my mind writing a book that was both physics and mathematics and that amalgamated all these ideas that were floating around. I had the feeling that we were at a point that was terminal in the sense to go any further we would need different techniques so I thought if I could summarise what the situation was right then it would remain in stone for some considerable time. I'd struggled to get started on this by myself and so I asked Ola what he thought about writing it together. I described the outline and how it would be half mathematics and half physics so we had complementary interests and he said yes and changed the subject but came back the next day and said we could give it a try. After Ola and I had started writing the book, maybe a couple of hundred pages, Thirring turned up in Marseilles. Someone told him the first evening that we were writing this book so first thing the next morning he came up to say I hear you are writing a book, I'm an editor at Springer, I'll accept it!

Later on the book project became a bit of a nightmare as I went to UNSW in Sydney in 1978 and Ola went back to Oslo. We had a very good exchange system by which I would write a section of a chapter and I'd send it to Ola and then he'd either say he didn't like it at all and he needed to rewrite it or that it was acceptable with minor corrections and whatever his decision it would come back to me under two weeks later. Now all the mail went through the international mail office in Sydney and they all went on strike for 3 months!

I was at UNSW from 1978-1982 and then moved to the Australian National University in 1982.

I never lost my love of cycling. During the time I was working at ANU I often took the 'long way' via Tuggeranong to ride from our house in Macquarie to the ANU and a similar length detour on the way home. You don't necessarily peak as a cyclist when you are young, I found I was faster at time trials at the age of 50 than I had been as an eighteen year old. In 2002 I was very pleased to win the World Masters Games 20k cycling race in the my age range. This was a time trial and I was the last person to start as I was the top seed and I completed the course in a little over 32 minutes. In some ways I consider this my greatest achievement!

\section*{Bratteli–Robinson 1974–2014}\label{D2}

Notes written by Derek after Ola Bratteli's death

Derek W. Robinson

Australian National University, Canberra

January 2017

My intention is to describe some details of the nearly forty years that I collaborated
with Ola. The talk will be largely historical rather than mathematical although it will
be difficult to avoid some mathematics. Since I am not a historian it will mainly be an
informal colloquial description of years gone by, not necessarily in the order they occurred.

I know little of formal history but I did once hear an after dinner speech by a distinguished
British historian who emphasized there were only three things to know about
history.

The first point about history is that ‘it is all behind us’, which is certainly true of my
collaboration with Ola. Truly and very regrettably. Initially neither of us ever envisaged
that we would form an enduring research partnership and I regret that it came to a premature
end two years ago. I also deeply regret losing so many other friends and colleagues in
the last two or three years, friends and colleagues who played a role in our collaboration:
Rudolf Haag, Daniel Kastler, Bill Arveson, Uffe Haagerup, John Roberts, Oscar Lanford,
Alan McIntosh.

The second point is that ‘history has no plot’ and that can certainly be said of the
collaboration with Ola. When we met we had quite different backgrounds, very different
personalities and totally different educations. None of these differences affected, however,
our ability to collaborate efficiently and it was always a great pleasure working with Ola.

I have forgotten the third all important historical rule. But maybe it was that history
is unreliable because it largely depends on the memories of the protagonists and these
memories are often coloured by subsequent reconstructions and interpretations. Certainly
the following account depends on my memory of distant events and the dates and details
should not be totally trusted.

I do not remember exactly when I first met Ola, it was sometime at the beginning of
1974. I nevertheless remember exactly where I met him. It was in the foyer, if you could
describe the bottom of a staircase as a foyer, in the then home of the CNRS Department of
Theoretical Physics in Marseille. This was a small building just off the Boulevarde Michelet
and not far from the famous Corbusier Unit´e d’Habitation. At that time I was Professor of
Physics at the Luminy Campus of the Universit\'{e} d’Aix-Marseille. This campus had been
founded shortly after `les \'{e}v\'{e}nements de mai' in 1968. In fact it only came into being in
1969 some six months after I had been appointed. Hence there were no offices initially,
in fact there were very few buildings, and that was why the theoretical physics group was
housed at the CNRS site.

Ola was introduced to me by Trond (Digernes) whom I had met earlier in 1971 in Los
Angeles. So that is probably a good point to start the story. In fact it was a good year
from many points of view. Ola was a graduate student in Oslo with Erling (Størmer) in
1971 and completed his now famous and often cited work on AF-algebras. It appeared in
the September 1972 edition of the Transactions of the American Mathematical Society. He
then continued his graduate studies at NYU. Trond at that time was a graduate student in
Los Angeles with Masimichi (Takesaki) working on operator algebras. I, on the other hand,
was entangled in several different areas of mathematics and physics, e.g. operator algebras,
quantum statistical mechanics, spin systems etc. During the 1971 trip to America when I
met Trond I wrote a paper on UHF-algebras, with Elliott Lieb, which was later to become
oft-cited in the physics literature, but that took another 35 years.

You might well wonder what two Norwegian graduate students of mathematics and
an English theoretical physicist were doing together in Marseille at that time in 1972 and
how did they come to have common interests in operator algebras. To explain that I have
to jump back a few more years. The starting point for our common interest in operator
algebras and theoretical physics began for me in the late 50s when I was a graduate student
in Oxford working on nuclear physics. I read a lot during that time and as a consequence
my personal interests turned toward quantum field theory. One paper I read, in the Journal
of the Danish Royal Society, which particularly fascinated me was by Rudolf Haag.
It considered representations of the algebra associated with fields satisfying the canonical
commutation relations. A key abstract idea was that the representation describing interacting
particles could not be unitarily equivalent to the representation describing free
particles. This clarified an earlier, more specific, result of van Hove. It was the bud that
bloomed slowly into algebraic quantum field theory. Although I tried to understand these
developments I had no opportunity to work in this area until I had completed my rather
mediocre thesis on the structure of deformed nuclei. That was in 1960. But I then had
the good fortune to spend two years as a NATO postdoctoral fellow at the Eidgen\"{o}ssiche
Technische Hochscule in Z\"{u}rich. In 1960–61 there was a small but active group headed by
Res Jost working on quantum field theory. David Ruelle also had a postdoctoral position
at the ETH during the same period and in 1961 his interest turned to obtaining precise
mathematical results in statistical mechanics. Simultaneously Rudolf Haag had also ventured
into this area and had given a precise analysis of the Bardeen–Cooper–Schrieffer
model of superconductivity based on algebraic representation theory. At that time there
no particular interest in algebraic methods in Z\"{u}rich although this changed in 1961–62
when Huzihiro Araki visited Z\"{u}rich and gave a lecture series on von Neumann algebras.
Araki completed his doctoral work with Rudolf Haag in 1960 at Princeton and his thesis
contained a detailed development of Haag’s 1955 theorem giving a precise link between
different interactions and distinct representations. The exposure to these various sources
was a strong influence on my subsequent work.

After the completion of my NATO fellowship in 1962 I had the very good fortune to
meet and work with Rudolf in Illinois. Quite coincidentally Araki also spent the year
1962–1963 in Rudolf’s group as did Daniel Kastler. Daniel had a sabbatical year from
his position as Professor of Theoretical Physics in Marseille. He had met Rudolf at the
1958 Varenna Summer School and was impressed with the idea of exploiting algebraic
methods to understand the local structure of quantum field theory. As a consequence
he had invested considerable time and energy reading the Murray–von Neumann papers
and learning all he could about the theory of $C^*$- and $W^*$-algebras. The fruits of his
labours were then harvested in the 1964 paper with Rudolf ‘An Algebraic Approach to
Quantum Field Theory’ which attracted a great deal of attention from both the physics
and mathematics communities. This paper had a significant impact which went beyond
its mathematical content since it reinvigorated Daniel’s efforts in developing mathematical
physics in Marseille. A development which led ten years later to Ola and I meeting at the
CNRS.

Let me at this point leap forward a bit in time, until 1968 when I moved to Marseille
for the beginning of a ten year stay. In 1965 I had visited for a couple of days and then
later that year I began a six month visit which led to Daniel campaigning for a permanent
position. In fact in 1968 it was not legally possible for a foreigner to be appointed to a
tenured position in a French University but this changed with the reforms following the
uprisings of 1968 and in September of that year I moved to Marseille. To be exact I moved
to Bandol, the small coastal village 45kms east of Marseille where Daniel and his family
lived. Subsequently my wife, Marion, and I built a house outside of Bandol on a hill
overlooking the mediterranean. It was an ideal setting for mathematics. By 1968 Daniel
had been successful in establishing his visitors program and there began to be a regular
flux of visitor from all parts. Erling was one of the earlier ones, in 1967 I believe. By good
fortune I was also visiting at that time from my then position at CERN. (I still had an
interest in particle physics.) Organizing a visitor’s program in the French system presented
many practical problems in those years, and possibly this has not changed. All paperwork
had to go through Paris and that led to unpredictable delays and errors. (When I was
in Marseille for six months in the 1965-66 period as a Professeur Associ\'{e} my salary did
not get paid until four and a half months had passed.) But as time went by we had an
increasing number of longer term visitors on sabbatical leave or travelling scholarships and
fellowships. Ola and Trond fell into the latter category with support from the beneficent
Norwegian government.

At the time I first met Trond after his arrival in Marseille in 1974 I was thinking
about derivations of $C^*$-algebras as these seemed significant for the description of physical
symmetries. The main problem was that the only general mathematical theorems for
derivations were, at that time, for bounded derivations. These were of little interest in
the context of symmetries since the corresponding generators were analogues of differential
operators and consequently were unbounded. Hence I was trying to understand what
were the general features of unbounded derivations that could be useful. I knew that
densely-defined symmetric operators on Hilbert space were automatically closeable and I
thought that there was possibly an algebraic analogue. Could it be that a densely-defined
symmetric derivation on a $C^*$-algebra is automatically closeable? I thought that this was
not unreasonable and asked Trond if he knew the answer. After assuring himself that it
was not obvious he said he would think about it. The following week I ran into him in the
foyer I mentioned earlier in the company of another person whom he introduced as Ola
Bratteli. Trond also immediately followed up this introduction by mentioning that he had
referred my question to Ola who had constructed a counterexample. This had me excited,
not a lot but a little, since I was expecting a positive rather than a negative example. Up
to this point Ola had said very little but then he broke into an explanation which I did
not fully understand, it certainly involved Cantor sets, of which I had scant knowledge,
and also seemed to be for abelian algebras. Somewhat later that afternoon he explained
things in a bit more detail but since I was primarily interested in non-abelian algebras,
even relatively simple ones such as UHF-algebras, I pressed him about more complicated
situations. That was the beginning of our collaboration.

Things went quickly. In a few weeks we had written most of our first paper on unbounded derivations of $C^*$-algebras. One of the themes of the paper was closeability and
we gave a couple of criteria and also examples of non-closeable derivations. The principal
example we gave was on a UHF-algebra generated by an increasing sequence of matrix
algebras on which the derivation was identically zero although the derivation itself was not
zero. This was constructed by extension of Ola’s original argument for abelian algebras.
Several years later Ola mentioned to me that he regretted that he had not written a more
comprehensive description of derivations on abelian algebras. It was certainly something
he understood well and which various other author wrote about some years later. I regret
that I did not encourage him more but I was focused on non-abelian situations and
quantum mechanical applications.

Our first paper also gave criteria for a derivation to generate a strongly continuous
group of $\null^*$-automorphisms. The latter were in large part adaptations of standard results
of semigroup theory. The special feature which, in hindsight, I realize we did not fully
appreciate was positivity. In fact semigroup theory which was widely viewed as a welldeveloped
and largely complete theory was severely deficient in respect to the analysis
of positivity properties. The one notable exception was a 1962 paper by Ralph Phillips
published in the Czechoslovakian Mathematical Journal which had apparently lapsed into
obscurity. This paper gave a very nice version of the Hille-Yosida theorem tailored to
positive semigroups. We only discovered this paper in 1980 and analyzed its implications
for C0-semigroups on $C^*$-algebras in a paper in Mathematica Scandinavica. Somewhat later
in the 1980s there were many other developments in the analysis of positivity properties.
In particular it was realized that a 1973 inequality of Kato for the Laplacian could be
adapted to give a criterion of positivity for quite general semigroups.

We were not the only people interested in unbounded derivations in the 1970s. Powers
and Sakai were working in this area but mainly on UHF algebras. Their belief was that each
symmetric derivation on such algebras was the limit of bounded derivations. Alternatively
stated the derivation was the asymptotic limit of inner derivations. This is a topic that
Akitaka will talk about after lunch so I will not go into further details. But the principal
method that was used, at least initially, was to analyze the functional properties of the
domain of the derivations. Our first paper also addressed this problem and we showed
that if $A = A^*$ was in the domain of a closed derivation $\delta$ on an abelian $C^*$-algebra
then $f(A)$ was also in the domain for each continuously differentiable function f on the
real line. In our second paper we gave a somewhat weaker statement that was valid for
general $C^*$-algebras with identity and was based on the observation that the identity was
automatically in the domain of $\delta$ and $\delta(\bdone)=0$. In the interim Bob Powers, who had gained
fame with the construction of a one-parameter family of non-isomorphic Type III factors in
his PhD thesis, had extended our abelian result to general $C^*$-algebras but unfortunately
his argument was incomplete. I was able to construct a counterexample to one vital step
in his proof. Two years later, on my first visit to Australia, I mentioned the problem to
Alan McIntosh and he gave a counterexample to the general statement. Fortunately Bob
forgave me pointing to the error of his ways and we later became good friends.

In the period 1974–1976 that Ola and I wrote our first three joint papers we were
in principle both in Marseille but in practise we were often travelling. I remember that
a significant part of our second paper was developed when I was in California visiting
Berkeley and Ola was in Marseille. It was during this period I met Bill Arveson and we
also became good friends. It was then that I met Tosio Kato for the first time. His book on
Perturbation Theory had been very influential for me. In particular I learnt almost all that
I then knew about semigroup theory from this book. In those days, of course, there were
no computers and no email so the only means of communication that Ola and I had was
ordinary mail. This certainly slowed things down but it was to become useful experience
for us when we later wrote the second volume of our book. At that point I was in Australia
and Ola was in Norway. But I will come to that a bit later. Another thing that slowed us
down in those early days was an accident that I had in which I squashed three vertebrae
in my back. This put me in hospital for a short period but kept me in bed, supine, for
several weeks. The hospital in Toulon wanted to fit me with a jacket that would keep my
back rigid but it turned out that I was too tall for all the jackets available. Therefore
they constructed a plaster cast which encased the whole upper part of my body. As it was
seriously warm in the South of France they left a large round aperture in the front of the
casing for cooling purposes. I resembled a frontloader washing machine with legs. I had to
wear this construction for five months and the weeks in bed were decidedly uncomfortable.
But as Autumn arrived I could move around and conduct business as usual. It was during
the period I was bedridden that Ola and I worked on our fourth paper. He often visited
Marion and I in the period I was confined to bed and although I was limited in movement
we found that was no impediment to doing mathematics.

By Autumn 1975 I was essentially recovered, my cast was removed and I could resume
a normal life. It was at this point I received an invitation which was to change my life
and, less directly, Ola’s life. The invitation was to lecture at a Summer Research School
organized by the Australian Mathematical Society in Adelaide. The prime mover behind
this invitation was Angas Hurst whom I had met at a similar summer meeting in 1961
in Yugoslavia. Since I had never been to Australia, although my parents flirted with the
idea of emigrating there in the years after the second world war, I was very excited by the
prospect of the trip. I was also apprehensive of sitting in an aeroplane for 24 hours so soon
after my accident. As it turned out the total journey from Marseille to Adelaide took 42
hours as I had to detour via London and Sydney. In Adelaide I lectured on various topics
including the work with Ola on derivations but mainly on the more physical aspects of my
work. I returned to Marseille in February 1976 excited by all the new experiences and also
seriously overweight partly as a consequence of my prolonged recuperation period but also
because of the excellent Australian food and wine. The latter was a particular surprise as
Australian wine was not at that time widely known in Europe and definitely not in France.
At that period I believe that I was probably heavier than Ola although the French food
and wine helped him forge ahead.

Two things happened after my return from Australia. First, full of enthusiasm for Australia,
Marion and I began to think of moving there at some point in the future and I wrote
to a couple of people about possibilities. This was foreseen as a longterm project. Secondly,
after my lectures in Adelaide Angas had begun to encourage me to write something
more extended about operator algebras and mathematical physics. In fact I had already
attempted that during a months stay in Groningen in 1972, a month during which Marion
and I simultaneously gave up smoking. It was a memorable month for the latter but not for
the former. I did fill three exercise books with draft material for a book but I realized that
I was not sufficiently prepared to write the mathematical background. Then around June
76 I suggested to Ola that we should write something together. The suggestion was a bit of
a surprise to him but within a few days he agreed. This was intended as a relatively short
term project. My idea was to write something about 3-400 pages with a couple of chapters
on mathematical background and a couple of chapters on applications to physics. Since the
major applications were at that point to models of statistical mechanics the book would
be Operator Algebras and Quantum Statistical Mechanics. We started work in September
76. But the best laid plans of mice and men oft go awry, to paraphrase Robbie Burns.

In October 76 the University of New SouthWales advertised a Chair in Pure Mathematics,
I applied, in February 77 I flew out to Sydney for an interview and shortly thereafter
I was offered the position. So the longterm project of moving to Australia became a very
short term project and we started to plan our move for January 78. During this period
Ola and I were working and writing and by September 77 we had the equivalent of 500
printed pages of material which exceeded somewhat our estimated length. That was the
good news. The bad news was that we were only half way through the planned material.
So the short term project turned out to be a long term project and the book changed from
one volume to two. It also meant that the second volume was largely written with Ola in
the Northern Hemisphere and me in the Southern Hemisphere.

Now I will try to explain how we planned, organized and wrote the book. Initially
we agreed on the outline as far as individual chapters were concerned. This outline was a
modification of the ideas I had earlier in Groningen. The final plan was for six chapters, one
a brief introduction on the background of the material to be covered, three on mathematical
topics and two on the more physical aspects. We started immediately with the second
chapter on the general theory of $C^*$-algebras and von Neumann algebras. We first made a
tentative sketch of the intended sections. Next we each took primary responsibility for one
or other of the sections. Then we would discuss the general presentation of the material in
each section. After these preliminaries we would begin to write drafts independently. For
example, I would draft the first section on general algebraic structure and Ola would draft
the second section on representations of algebras. Then we would swap the drafts and each
edit the others work. This process would be repeated until we were each satisfied with the
outcome. I am not sure whether this is a standard procedure with coauthors of books
but it worked well for us. The editing was not a superficial process since we often had
different notions of the relative significance of the material and the emphasis to be given
to various statements and results but this would be ironed out in the various exchanges.
At times my first draft would be completely changed by Ola and vice versa. Somehow
the process always reached equilibrium after a reasonably short time, with one exception.
This procedure also had various advantages. It naturally introduced a uniformity of style.
It also gave a fairly foolproof method of avoiding error, although we were not completely
successful in that respect.

There were two style decisions that we agreed on at the beginning. The first was to
eschew footnotes. This artifice is used by many authors as a means of citing the relevant
literature. We, however, decided that we would not give references in the main text but
add notes on the background and references in sections at the end of each chapter. This
also allowed us to introduce some digressions from the main themes. The problem with
this method is that one can easily offend colleagues by false attribution of priority or by
oversight of their contribution. I am not sure whether we did offend anyone but we never
had any complaints so I think our ‘history’ was probably reliable.

The section which caused the most difficulty was the section on Tomita–Takesaki theory.
This was a very opaque theory at the time but Ola made a brave attempt to describe it
in a transparent way. Unfortunately I did not find the end result at all transparent so I
completely rewrote it. Then Ola was dissatisfied with my presentation and rewrote it once
more. Finally we had seven drafts before we came to an agreement. We were both pleased
a couple of years later to hear that Alain Connes was recommending our description as the
best introduction to the material.

The collaborative exchange process we developed for the chapter on the general algebraic
theory was continued for all the other chapters. It meant that all parts of the book,
with a couple of exceptions were really written by both of us. It also meant that we were
not guilty of just transcribing proofs of other authors. But we did run into some problems
of presentation even where least expected. For example our third chapter was on the theory
of one-parameter groups and semigroups of operators. This theory was well developed
for strongly continuous groups and was immediately applicable to the $C^*$-algebra theory.
But for $C^*$-algebras one needed to consider weak$^*$-continuous semigroups. The change of
topology affected several of the standard results and we had to sort this out. Ola then
made the suggestion that we should develop the theory in the language of dual topologies
and this gave a unity to the description. It is only in the last five years that I have seen
another genuine attempt at unification. In addition Ola suggested we should consider Jordan
algebras and positive maps. I think he had learnt about such things through his days
as a student of Erling. He was well versed in this theory but it was a topic I learned in the
writing.

A different type of problem occurred with Chapter 4 on Decomposition Theory. The
theory of decomposition of representations was well established but the theory of decomposition
of states of operator algebras was much more recent. There had been lots of
developments in the late 60s and early 70s and I had been heavily involved with several
of them. The literature was a melange of Choquet theory, operator theory, functional
analysis, ergodic theory etc. involving central, subcentral, extremal and ergodic decompositions.
We had to try and unify results by many authors who had used disparate notations,
definitions and techniques. Both Ola and I were very satisfied with the end result.

The theory of decomposition of representations was a different kettle of fish. This
originated with von Neumann in the late 40s and had undergone very little change. I was
not sure whether it was worthwhile including a description as it was already well covered
in other books. But Ola bravely volunteered to write a version and although I read it and
made a few minor changes it is all his own work.

We had finished writing Chapters 2, 3 and 4 by August 77 and then Ola went on
holiday. So in his absence I also did a bit of privateering and wrote Chapter 1. When Ola
returned and I gave him my draft manuscript he let out an enormous sigh of relief and
admitted that it was the one thing he had not looked forward to writing. He did suggest
a few small changes but other than that I can claim it as a personal contribution. These
were the only two deviations from our combined, unified, writing method.

We realized that summer that we had begun to write the book with no specific plans
for its publication so we began to discuss possible publishers. But any potential difficulties
were overcome by a visit of Walter Thirring to Marseille. He apparently heard about our
project and immediately came to my office and suggested he could have it published in
the Springer series, Texts and Monographs in Physics, of which he was an editor. Since he
had not seen any of the manuscript this was a very flattering offer that provided an easy
solution to the publication problem. In fact we subsequently had a couple of difficulties
with the publishing house which left us doubting whether we had made the correct choice.
Anyway we did proceed and submitted the first four chapters as the first volume of the
book.

In September 77 we began to see the problems ahead. I was due to leave Marseille at
Christmas and so that only left us a bit over 3 months to write the second volume. On top
of that I had to teach full time for two of those months. (One of the reforms introduced
in Luminy after the ‘revolution’ of 1968 was to introduce block teaching. A typical course
would be taught for sixteen hours a week over two months instead of four hours a week
over eight months.) But we decided to do as much as possible in the time remaining.

Chapter 5, the first chapter of Volume 2, was a mixture of well-established material,
i.e. the free bose gas and the free Fermi gas, and new material only developed in recent
years, i.e. the structure of KMS states. Somewhat miraculously we managed to write the
major part of this chapter before I was due to leave. One thing that helped was the
arrival of Akitaka Kishimoto and he was an enormous help with the analysis of stability
properties of KMS states. This led to a separate research paper the main results of which
were integrated into the book. In retrospect I am amazed that we accomplished so much
in such a short time but we nevertheless were far from finishing the whole book. Parts of
Chapter 5 and all of Chapter 6 remained.

My recollections of the first half of 1978 are rather blurred. Marion and I, together
with our two young daughters, arrived in Sydney late in January and the first semester
of the Australian academic year began toward the end of February. I think the book
writing came to a halt. But I also believe the proofs of Volume 1 arrived in this period
and they provided a considerable amount of work. They also provided one major surprise.
Springer had unilaterally decided that all proofs would be in a different smaller typeface
than the general text. Unfortunately they had never asked us to mark the ends of proofs.
Therefore the decision on the points at which the typeface should change had been made
by some member of the Springer editorial staff. They clearly did not understand the text.
There were the inevitable errors to correct in the proofs but large expanses of typeface had
to be changed. This was still in the era that the typesetting was done in the traditional
Gutenberg manner so this involved a considerable amount of extra work and cost. Springer
expected us to pay for it. This was the first point that I had doubts about our choice of
publisher. I remember corresponding with Ola about the problem and he agreed that I
should try and negotiate with Springer. I was not in the mood for any form of friendly
negotiation and my letters to Springer and the entire editorial board of the Texts and
Monographs series reflected my vehemence. Fortunately Springer accepted that the error
was on their part and waived the costs. There was one other unfortunate consequence of
the typeface incident. Reading so much small type in such a short period caused my eyes
to deteriorate and for the first time in my life I needed glasses.

The next time I saw Ola was in the Australian midyear break, which would have been
late June, early July 78. I returned to Marseille and our collaboration started up again. I
recently checked that the paper with Akitaka and Ola on stability and the KMS condition
was submitted for publication at the beginning of March 78 so I assume the version that
appears in Chapter 5 had been written earlier or that we wrote it that month. Akitaka was
still in Marseille and the three of us wrote a paper on the ground states of quantum spin
systems. Ground states were a topic of Chapter 5 and a section of Chapter 6 dealt with
ground states of spin systems. So I infer that we were completing Chapter 5 and thinking
ahead for Chapter 6. Certainly after my return to Australia we began to seriously write
about spin systems. About the time I returned to Australia Ola and Akitaka also left
Marseille. Ola returned to Oslo and Akitaka moved to Ottawa.

The next eleven months Ola and I continued our collaboration regularly by mail. Remember
that there were no convenient personal computers at that time and certainly no
email. Although Steve Jobs and Steve Wozniak were producing the first prototype of the
Apple computer in a garage in Colorado and the revolution in communication was on the
horizon. Unfortunately it was a bit too late to help us with the task of completing Chapter
6. Fortunately there was a a reliable airmail system between Australia and Norway.
Typically one of us would write a section of manuscript mail it to the other who would
edit it and mail it back. The whole exchange took two weeks. It was not ideal and at
times it was very frustrating as it broke up the continuity. Nevertheless we were able to
proceed with this system and we adapted to the routine. Then, however, there was a
disaster. International airmail passing through Sydney was handled at a postal exchange
in the centre of the city and at the beginning of 1979 the workers all went on strike. I
have no memory of their grievances, justified or not, but the strike lasted two months and
during this period there was simply no overseas mail. Again our collaboration ground to
a halt. But life went on.

After the end of the Sydney postal strike our collaboration rebooted and we made
serious progress on the final chapter of the book. Then in the midyear break we both
returned to Marseille and began the final work. In fact we set up in Bandol in one of
the Katikias apartments. We had four weeks to complete the book and it took us three
and a half. Although the book was in some sense completed there remained a great deal
more drudgery. There was no TEX in those days and typing mathematics depended on
the now archaic typewriters with interchangeable balls which allowed for a wide variety
of mathematical notation. It was impractical to prepare any kind of usuable draft by
this method and all manuscripts had to be written by hand and then typed by an expert
typist, an expert with a great deal of patience. I believe the text we wrote that summer
was eventually typed back in Sydney. Anyway later in the (northern) Autumn of 79 the
completed manuscript was submitted to Springer and the second volume finally appeared
in 1981. So the total operation of writing and publishing the 1000 page two volume book
took about three and a half years.

There was one topic in Chapter 6 that I referred to earlier in the talk and which
provides a salutary warning to the research managers who rely on citation indices to assess
performance. Toward the end of Section 6.2.1 we derive and discuss a property of finite
speed of propagation for a large class of quantum spin systems. This material was an
expansion of a paper I wrote with Elliott Lieb in 1972, the week after I met Trond for
the first time, and a subsequent paper of mine in 1976. These papers went down like the
proverbial lead balloon. Other than self-citations and a mention by Park they were totally
ignored until 2006. Since then the paper with Lieb has been cited approximately once a
week and it even has its own Wikipedia page. The explanation–the development of the
theory of Quantum Computing. So don’t despair if your work is not cited immediately.

Although the manuscript of the book was finished in 1979 that did not slow down the
collaboration with Ola and Akitaka. They both visited Sydney for extended periods and we
wrote several papers on the properties of positivity preserving semigroups. This topic was
to receive a great deal of attention during the 80s and it is still a fruitful area of research.
Of course the generators of positive contractive semigroups are determined by Dirichlet
forms and the theory of the latter expanded rapidly in the same period and continues to
expand. From 1979 onward Ola was a regular visitor the Australia. He realized that by
timing the visits correctly he could ensure that it was always summer.

In late 1981 Marion and I were on the move again but this time it was only a short
(300 kilometers) hop to Canberra where we have remained ever since. In the first years
in Canberra we had many visitors working in diverse areas. Ola, Akitaka and I kept on
working on aspects of operator algebras and somewhat less on continuous semigroups. By
1985 we had recovered from writing the book, which was still selling regularly after the
first major sales to university libraries. But then Ola and I independently noticed that
the sales of Volume 1 ceased. We did not understand why until Bill Arveson’s house in
the hills overlooking Berkeley burnt down. You might well think that these events were
unrelated but Bill had lots of his mathematics books in the house and when he tried to buy
a replacement for our book he was told that Volume 1 was not available. He then wrote to
Ola asking for his assistance in locating a copy, Ola informed me and we decided to write
to Springer asking why the volume was out of print. The explanation we received again
gave us doubts about our rapid choice of publisher. Apparently Springer decided to close
their New York office and also decided to destroy all copies of books in underperforming
series which were stockpiled in New York. Our book which was the best performing in the
Physics series was a victim of the group action.

Ola and I were both horrified by this revelation and wrote that we would like Volume 1
to be reprinted. Springer then made a counter offer that they would publish a revised and
updated version as a second edition. So we had to start again. In some ways it was a good
thing. It gave us an opportunity to correct two small errors and also to add some discussion
of new developments. The latter was not difficult as not much had changed in six years
that had intervened. The subject matter was in a stable state. It was also not difficult
since we were both still working in the same general area. But things were somewhat
different in 1995 when Springer asked us to produce a second edition of Volume 2.

By the mid 90s I was principally working on elliptic operators and subelliptic operators
and Ola, Palle Jørgensen, Charles Batty also participated on occasions. Our interests
were far from models of quantum statistical mechanics. But we were fortunate in three
respects. First Ola was aware of the work of Bost and Connes on KMS systems and
this allowed us to include a description of their theory. Secondly I had followed some of
the developments in the theory of phase transitions in spin systems so that allowed us to
improve our original description of these topics and move the original emphasis away from
classical systems to quantum systems. Thirdly, MathSciNet arrived on the scene. As Ola
was arriving in Canberra for a months stay to revise Volume 2 the university obtained free
access for a month to MathSciNet. I did not realize the significance of this development
but I mentioned it to Ola and he immediately started searching for all work that had been
published in the preceding 15 years which might be of interest. We were amazed how
helpful the new technology could be. Although we started to prepare the second edition of
Volume 2 with great trepidation we were finally satisfied with our efforts. A few years later
when a new editor with Springer suggested we write a third volume we did not hesitate in
refusing.

At this point I will break off from the temporal description of our collaboration. Ola
continued to visit Australia in his attempt to make life a continuous summer and he often
visited our beach house. Canberra is about a two hour drive from the Pacific coast and
the nearest coastal town, Bateman’s Bay, lies at the mouth of the Shoalhaven River. The
latter is renowned for its oysters and they were one of Ola’s main delights. He became
a well-valued customer at the main oyster outlet in town. He also was happy to spend
days on the beaches which by any European standard were almost deserted. He enjoyed
swimming and would cover much longer distances than I ever fancied. After we sold the
beach house we rarely visited the coast but Ola would drive down for the day to have fish
for lunch and to bring back a bag of oysters. The last time I saw him was at Oslo Airport
and I suggested that he might be able to visit Australia and Bateman’s Bay once again
but he shook his head wryly. Unfortunately it was never possible.

%
%



\end{document}